\documentclass[a4paper,10pt]{article}

\usepackage{amsmath,amssymb,amsthm,xypic,epsfig}
\input xy
\xyoption{all}

\let\al\alpha

 \let\Dl\Delta
\let\epe\epsilon \let\eps\varepsilon \let\epsilon\eps

\let\la\lambda 
\let\om\omega 
 \let\phi\varphi

\newcommand{\Z}{{\Bbb Z}}

\newcommand{\R}{{\Bbb R}}
\newcommand{\C}{{\Bbb C}}

\newcommand{\Ref}[1]{{$($\ref{#1}$)$}}
\newcommand{\bean}{\begin{eqnarray}}
\newcommand{\eean}{\end{eqnarray}}
\newcommand{\be}{\begin{displaymath}}
\newcommand{\ee}{\end{displaymath}}
\newcommand{\bea}{\begin{eqnarray*}}
\newcommand{\eea}{\end{eqnarray*}}
\newcommand{\g}{{{\frak g}\,}}

\newcommand{\n}{{{\frak n}}}
\newcommand{\h}{{{\frak h\,}}}
\newcommand{\Id}{{\operatorname{Id}}}

\newcommand{\T}{\!\otimes\!}

\newcommand{\vs}{\vspace{.5\baselineskip}}

\newtheorem%
{thm}{Theorem}
\newtheorem%
{proposition}[thm]{Proposition}
\newtheorem%
{lemma}[thm]{Lemma}
\newtheorem%
{lemmadef}[thm]{Lemma-Definition}
\newtheorem%
{corollary}[thm]{Corollary}
\newtheorem%
{conjecture}[thm]{Conjecture}
\newcommand{\End}{{\operatorname{End\,}}}
\newcommand{\Hom}{{\operatorname{Hom\,}}}

\newcommand{\wo}{{\widetilde \otimes}}
\newcommand{\bo}{{\bar {\otimes}}}

\newcommand{\Ti}{{\cal {T}}}

\newcommand{\Tb}{{\tilde { \otimes}\,}}

\def\vsk#1>{\vskip#1\baselineskip}
\def\ftext#1{{\let\thefootnote\relax\footnotetext{\vsk-.8>\noindent #1}}}

\def\g{\mathfrak{g}}
\def\h{\mathfrak{h}}
\def\n{\mathfrak{n}}
\def\b{\mathfrak{b}}
\def\C{\mathbb{C}}
\def\R{\mathbb{R}}
\def\Z{\mathbb{Z}}
\def\qed{$\hfill \blacksquare$}
\def\l{\mathfrak{l}}
\def\Ro{\mathcal{R}_0}

\newtheorem{theo}{Theorem}[section]
\newtheorem{prop}{Proposition}[section]
\newtheorem{lem}{Lemma}[section]
\newtheorem{cor}{Corollary}[section]
\numberwithin{equation}{section}

\begin{document}
\title{\huge{{Lectures on 
the dynamical Yang-Baxter equations}}
}
\author{Pavel Etingof and Olivier Schiffmann}
\date{}
\maketitle

\section{Introduction}

This paper arose from a minicourse given by the first author 
at MIT in the Spring of 1999, when the second author extended 
and improved his lecture notes of this minicourse. 
It contains a systematic and elementary 
introduction to a new area of the theory of quantum groups --
the theory of the classical and quantum dynamical Yang-Baxter equations. 

The quantum dynamical Yang-Baxter equation is a generalization 
of the ordinary quantum Yang-Baxter equation. It first appeared
in physical literature in the work of Gervais and Neveu \cite{GN}, 
and was first considered from a mathematical viewpoint by Felder
\cite{F}, who attached to every solution of this equation a quantum group, 
and an interesting system of difference equations, -
the quantum Knizhnik-Zamolodchikov-Bernard (qKZB) equation.
Felder also considered the classical analogue of the 
quantum dynamical Yang-Baxter equation -- the classical 
dynamical Yang-Baxter equation. Since then, this theory 
was systematically developed in many papers, some of which are 
listed below. 
By now, the theory of the classical and quantum dynamical Yang-Baxter 
 equations and their solutions has
many applications, in particular to integrable systems and representation 
theory. To discuss this theory and some of its applications 
is the goal of this paper.   
 
The structure of the paper is as follows. 

In Section 2 we consider the exchange construction, which is a natural 
construction in classical representation theory that leads one to discover 
the quantum dynamical Yang-Baxter equation and interesting 
solutions of this equation (dynamical R-matrices). 
In this section we define the main objects 
of the paper -- the fusion and exchange matrices for Lie algebras and quantum 
groups, and compute them for the Lie 
algebra $sl_2$ and quantum group $U_q(sl_2)$.   

In Section 3 we define the quantum dynamical Yang-Baxter equation,
and see that the exchange matrices are solutions of this equation. 
We also study the quasiclassical limit of the quantum dynamical
Yang-Baxter  equation -- the classical dynamical Yang-Baxter
equation.  We conjecture that any solution of this equation can
be quantized.  We compute classical limits of exchange
matrices,  which provides interesting examples of solutions 
of the classical dynamical Yang-Baxter equation, which we call 
basic solutions. 

In Section 4 we give a classification of solutions of the classical 
dynamical Yang-Baxter equation for simple Lie algebras defined on a Cartan
subalgebra, satisfying the unitarity condition. 
The result is, roughly, that all such solutions 
can be obtained from the basic solutions.  

In Section 5 we discuss the geometric interpretation 
of solutions of the classical dynamical Yang-Baxter equation, which 
generalizes Drinfeld's geometric interpretation of solutions of the classical 
Yang-Baxter equation via Poisson-Lie groups. This interpretation is 
in terms of Poisson-Lie groupoids introduced by Weinstein. 

In Section 6 we give a classification of solutions 
of the quantum dynamical Yang-Baxter equation for the vector representation 
of $gl_N$, satisfying the Hecke condition. As in the classical
case, the result states that  all such solutions can be obtained
from the basic solutions which arise from the exchange
construction. 

In Section 7 we discuss the "noncommutative geometric" interpretation 
of solutions of the quantum dynamical Yang-Baxter equation, which 
generalizes the interpretation of solutions of the quantum 
Yang-Baxter equation via quantum groups. This interpretation is 
in terms of quantum groupoids (or, more precisely, $H$-Hopf
algebroids). 
 
In Section 8 we give a defining equation satisfied by the 
universal fusion matrix --
the Arnaudon-Buffenoir-Ragoucy-Roche (ABRR) equation, and prove it 
in the Lie algebra case. We give applications 
of this equation to computing the quasiclassical limit of the fusion
matrix, and to computation of the fusion matrix itself for $sl_2$.    
  
In Section 9 we discuss the connection of solutions
of the quantum dynamical Yang-Baxter equation 
to integrable systems and special functions, in particular to 
Macdonald's theory. Namely, we consider weighted traces of intertwining 
operators between representations of quantum groups, and 
give difference equations for them which in a special case 
reduce to Macdonald-Ruijsenaars difference equations. 

Appendix A contains the classification of solutions 
of the classical dynamical Yang-Baxter equation for simple Lie algebras
defined on subspaces of the Cartan subalgebra. 

Appendix B contains a proof of the ABRR equation in the quantum case. 

In Appendix C we give a characterization of the set of solutions of the
\textit{dynamical 2-cocycle equation} (which is satisfied by the fusion
matrices) in terms of a purely algebraic equation.

Appendix D makes the link between fusion matrices and Shapovalov forms on
Verma modules over Lie algebras or quantum groups. 

At the end we review some of the existing literature that is relevant to the 
theory of the dynamical Yang-Baxter equations. 
 
To keep these lectures within bounds, we do not discuss dynamical Yang-Baxter 
equations with spectral parameter. These equations are related 
to affine Lie algebras and quantum affine algebras just like 
the equations without spectral parameter are related to finite
dimensional Lie  algebras and quantum groups. Most of the
definitions and results  of these lectures can be carried over to
this case, which gives  rise to a more interesting but also more
complicated theory  than the theory described here. A serious 
discussion of this theory would require 
a separate course of lectures. 

Detailed proofs of several statements made in these notes can be
found on T. Koornwinder's webpage at :
$$http://turing.wins.uva.nl/~thk/recentpapers/comment.html$$

{\bf Acknowledgements.} We thank the participants of the 
minicourse at MIT and of the "Quantum groups" conference in Durham 
(July 1999) for interesting remarks and discussions. We are
grateful to Ph. Roche and A. Varchenko for many useful
conversations, in particular for suggesting to discuss 
the precise relationship between the fusion matrix and the
Shapovalov form (Appendix D).
We thank IHES and Harvard University for hospitality.
The work of P.E. was partially supported by the NSF grant DMS-9700477,
and was partly done while he was employed by the Clay 
Mathematics Institute as a CMI Prize Fellow. We are indebted
to Tom Koornwinder for useful comments on this paper and for pointing out
several misprints in an earlier version of these notes (see the above
mentioned webpage); these misprints have been corrected in 
the present version.

\section{Intertwining operators, fusion and exchange matrices.}
\paragraph{2.1. The exchange construction.} We start by giving a
simple and natural  construction in classical representation
theory which leads to discovery of the quantum dynamical
Yang-Baxter equation. 

Let $\g$ be a simple complex Lie algebra, $\h\subset\g$ a Cartan
subalgebra and $\Delta \subset \h^*$ the associated root system. Let
$\Pi$ be a set of simple roots, $\Delta^+ \subset \Delta$ the
associated system of positive roots. Let $\g=\n_- \oplus \h \oplus
\n_+$ be the corresponding polarization of $\g$ and let $\g_\alpha$ be
the root subspaces of $\g$. Let $\langle\;,\;\rangle$ be the
nondegenerate invariant symmetric form on $\g$ normalized by the
condition $\langle \alpha,\alpha \rangle=2$ for long roots. Finally,
for each $\alpha \in \Delta$, choose some $e_\alpha \in \g_\alpha$ in
such a way that $\langle e_\alpha,e_{-\alpha} \rangle=1$.\\ \hbox
to1em{\hfill}For $\lambda \in \h^*$, let $\C_\lambda$ be the
one-dimensional $(\h\oplus \n_+)$-module such that $\C_\lambda=\C
x_\lambda$ with $h.x_\lambda=\lambda(h)x_\lambda$ for $h \in \h$ and
$\n_+.x_\lambda=0$. The Verma module of highest weight $\lambda$ is
the induced module
$$M_\lambda=\mathrm{Ind}_{\h\oplus\n_+}^\g\;\C_\lambda.$$

\paragraph{}Notice that $M_\lambda$ is a free $U(\n_-)$-module and
can be identified with $U(\n_-)$ as a linear space by the map $U(\n_-)
\stackrel{\sim}{\to} M_\lambda$, $u \mapsto u.x_\lambda$.\\ \hbox
to1em{\hfill}Define a partial order on $\h^*$ by putting $\mu<\nu$ if
there exist $\alpha_1,\ldots \alpha_r \in \Delta^+$, $r>0$, such
that
$\nu=\mu+\alpha_1+\ldots +\alpha_r$. Let $M_\lambda=\bigoplus_{\mu\leq
\lambda} M_\lambda[\mu]$ denote the decomposition of $M_\lambda$ into
weight subspaces.

The following proposition is standard.
\begin{prop}The module $M_\lambda$ is irreducible 
for generic values of $\lambda$.
\end{prop}

\paragraph{}Define also the dual Verma module $M^*_\lambda$ to be the
graded dual vector space $\bigoplus_\mu M_\lambda [\mu]^*$ equipped
with the following $\g$-action: $$(a.u)(v)=-u(a.v)\;\forall a \in
\g,\;u \in M_\lambda^*,\;v\in M_\lambda.$$ Let $x_\lambda^*$ be the
lowest weight vector of $M_\lambda^*$ satisfying $\langle
x_\lambda,x_\lambda^*\rangle=1$.  

\paragraph{}Now let $V$ be a
finite-dimensional $\g$-module.  Let $V=\bigoplus_{\nu \in
\h^*}V[\nu]$ be its decomposition into weight subspaces. Let
$\lambda,\mu \in \h^*$ and let us consider $\g$-module intertwining
operators $$\Phi:\;M_\lambda \to M_\mu \otimes V.$$ If $\Phi$ is such
an intertwining operator, define its "expectation value'' by $$\langle
\Phi \rangle=\langle \Phi.x_\lambda,x_{\mu}^*\rangle \in
V[\lambda-\mu].$$ 

\paragraph{Remark.} This definition is similar to the
notion of expectation value in quantum field theory.  

\begin{prop} Let
$M_\mu$ be irreducible. Then the map $$\mathrm{Hom}_\g
(M_\lambda,M_\mu \otimes V) \to V[\lambda-\mu],\; \Phi \mapsto \langle
\Phi \rangle$$ is an isomorphism.  
\end{prop} 

\noindent
\textit{Proof.} By Frobenius reciprocity, we have $$\mathrm{Hom}_\g
(M_\lambda,M_\mu \otimes
V)=\mathrm{Hom}_{\h\oplus\n_+}(\C_\lambda,M_\mu \otimes
V)=\mathrm{Hom}_{\h \oplus \n_+}(\C_\lambda \otimes M_\mu^*,V).$$
Moreover, since $M_\mu$ is irreducible, we have
$M_\mu^*=\mathrm{Ind}_{\h}^{\h \oplus \n_+} \C_{-\mu}$ as an $\h
\oplus \n_+$-module. In particular, $$\mathrm{Hom}_{\h \oplus
\n_+}(\C_\lambda \otimes M_\mu^*,V)=\mathrm{Hom}_\h (\C_\lambda
\otimes \C_{-\mu},V)=V[\lambda-\mu].$$ \qed 

\paragraph{}This
proposition can be reformulated as follows: for any $v \in
V[\lambda-\mu]$ there exists a unique intertwining operator
$\Phi^v_\lambda:\;M_\lambda \to M_\mu \otimes V$ such that
$$\Phi^v_\lambda(x_\lambda) \in x_\mu \otimes v + \bigoplus_{\nu <\mu}
M_\mu[\nu] \otimes V.$$ Notice that $\Phi^v_\lambda$ (for fixed $v$)
is defined only for generic values of $\lambda$. Identifying the Verma
modules $M_\lambda$ and $M_\mu$ with $U(\n_-)$, we can view
$\Phi^v_\lambda$ as a linear map $U(\n_-) \to U(\n_-) \otimes V$. It
is easy to see that the coefficients of this map (in any basis) are
rational functions of $\lambda$.  

\paragraph{}We would now like to
consider the "algebra" of such intertwining operators. Let us denote
by $\mathrm{wt}(u)\in\h^*$ the weight of any homogeneous vector $u$ in
a $\g$-module. Let $V,W$ be two finite-dimensional $\g$-modules, and
let $v \in V,\;w\in W$ be two homogeneous vectors. Let $\lambda \in
\h^*$ and consider the composition
$$ \Phi^{w,v}_\lambda:\;M_\lambda
\stackrel{\Phi^v_\lambda}{\longrightarrow} M_{\lambda-\mathrm{wt}(v)}
\otimes V
\stackrel{\Phi^w_{\lambda-\mathrm{wt}(v)}}{\longrightarrow}
M_{\lambda-\mathrm{wt}(v)-\mathrm{wt}(w)} \otimes W
\otimes V.$$ (Here and below we abuse notations and write
$\Phi$ instead of $\Phi \otimes 1$). Then
$\Phi^{w,v}_\lambda \in
\mathrm{Hom}_\g(M_\lambda,M_{\lambda-\mathrm{wt}(v)-\mathrm{wt}(w)}
\otimes W \otimes V)$. Hence by Proposition 2.2, for generic $\lambda$
there exists a unique element $u \in W \otimes
V[\mathrm{wt}(v)+\mathrm{wt}(w)]$ such that
$\Phi^u_\lambda=\Phi^{w,v}_\lambda$. It is clear that the
assignment
$(w,v) \mapsto u$ is bilinear, and defines an $\h$-linear map
\begin{align*}
J_{WV}(\lambda):\; W \otimes V &\to W \otimes
V,\\
w \otimes v &\mapsto \langle \Phi_\lambda^{w,v} \rangle
\end{align*}
\paragraph{Definition.} We call the operator
$J_{WV}(\lambda)$ the
\textit{fusion matrix} of $V$ and $W$.  

\paragraph{}We will now list
some fundamental properties of fusion matrices.  First let us
introduce an important piece of notation to be used throughout this
text. If $A_1,\ldots A_r$ are semisimple $\h$-modules and $F(\lambda):
A_1 \otimes \ldots \otimes A_r \to A_1 \otimes \ldots \otimes A_r$ is a
linear operator depending on $\lambda \in \h^*$ then, for any
homogeneous $a_1,\ldots a_r$ we set $$F(\lambda-h^{(i)})(a_1 \otimes
\ldots \otimes a_r):=F(\lambda-\mathrm{wt}(a_i)) (a_1 \otimes \ldots
\otimes a_r).$$ 

\begin{prop} Let $V,W$ be finite-dimensional
$\g$-modules. Then 
\begin{enumerate} \item $J_{WV}(\lambda)$ is a
rational function of $\lambda$.  \item $J_{WV}(\lambda)$ is strictly
lower triangular, i.e. $J=1+N$ where $$ N(W[\nu] \otimes V[\mu])
\subset \bigoplus_{\tau<\nu, \mu<\sigma} W[\tau]\otimes V[\sigma].$$
In particular, $J_{WV}(\lambda)$ is invertible.  \item Let $U,V,W$ be
finite-dimensional $\g$-modules. Then the fusion matrices satisfy the
following \text{dynamical 2-cocycle condition}:
 $$\qquad \qquad \qquad J_{U \otimes
W,V}(\lambda)(J_{UW}(\lambda-h^{(3)}) \otimes 1)=J_{U,W \otimes
V}(\lambda) (1 \otimes J_{WV}(\lambda)).  $$ on $U\otimes W\otimes V$.
\end{enumerate} 
\end{prop} 
\noindent \textit{Proof.} Statements 1. and
2. follow from the definitions and from the fact that the intertwining
operators $\Phi^v_\lambda$ are rational functions of $\lambda$. To
prove statement 3., let $u\in U, \;v \in V,\; w \in W$ be homogeneous
elements and consider the composition 
\begin{equation*}
\begin{split}
 M_\lambda
\stackrel{\Phi^v_\lambda}{\longrightarrow} M_{\lambda-\mathrm{wt}(v)}
\otimes V \stackrel{\Phi^w_{\lambda-\mathrm{wt}(v)}}{\longrightarrow}
&M_{\lambda-\mathrm{wt}(v)-\mathrm{wt}(w)} \otimes W \otimes V\\
&\stackrel{\Phi^u_{\lambda-\mathrm{wt}(v)-\mathrm{wt}(w)}}{\longrightarrow}
M_{\lambda-\mathrm{wt}(u)-\mathrm{wt}(v)-\mathrm{wt}(w)} \otimes U
\otimes W \otimes V.
\end{split}
\end{equation*}
The dynamical 2-cocycle condition follows from the
associativity relation
$$\Phi^u_{\lambda-\mathrm{wt}(v)-\mathrm{wt}(w)}\circ
(\Phi^w_{\lambda-\mathrm{wt}(v)}\circ\Phi^v_\lambda)=(\Phi^u_{\lambda-
\mathrm{wt}(v)-\mathrm{wt}(w)}\circ
\Phi^w_{\lambda-\mathrm{wt}(v)})\circ\Phi^v_\lambda$$ and from the
definition of the fusion matrices. \qed 
\paragraph{}The fusion
matrices can be viewed as the structure constants for multiplication
in the "algebra" of intertwining operators. We now turn to the
structure constants for "commutation relations". Let $V,W$ be two
finite-dimensional $\g$-modules. Let us define
$$R_{VW}(\lambda)=J_{VW}(\lambda)^{-1}J_{WV}^{21}(\lambda) \in
\mathrm{Hom}_\h (V \otimes W, V \otimes W),$$ where $J^{21}=PJP$ with
$P(x \otimes y)=y \otimes x$. The above definition can be rephrased in
terms of intertwining operators as follows: $R_{VW}(\lambda) (v
\otimes w)=\sum_i v_i \otimes w_i$ where $\Phi^{w,v}_\lambda=P\sum_i
\Phi_\lambda ^{v_i,w_i}$.

\paragraph{Definition.} The operator $R_{VW}(\lambda)$ is called the
\textit{exchange matrix} of $V$ and $W$.

\begin{prop} Let $U,V,W$ be three finite-dimensional
$\g$-modules. Then the exchange matrices satisfy the following
relation 
\begin{equation}\label{E:1}
R_{VW}(\lambda-h^{(3)})R_{VU}(\lambda)R_{WU}(\lambda-h^{(1)})=R_{WU}(\lambda)R_{VU}(\lambda-h^{(2)})R_{VW}(\lambda)
\end{equation}
 in the algebra $\mathrm{Hom}_\h(V \otimes W \otimes U, V \otimes W \otimes U)$.
\end{prop}
\noindent
\textit{Proof.} Let $u\in U,\; v\in V,\; w\in W$ be homogeneous elements and, as in Proposition 2.3, consider the
composition $\Phi_\lambda^{u,w,v}=\Phi^u_{\lambda-\mathrm{wt}(v)-\mathrm{wt}(w)}\circ
\Phi^w_{\lambda-\mathrm{wt}(v)}\circ\Phi^v_\lambda$. The proof of relation (\ref{E:1}) is obtained by
rewriting $\Phi^{u,w,v}_\lambda$ as  $\sum \sigma \Phi^{v_i,w_i,u_i}_\lambda$ where $\sigma: U \otimes W
\otimes V \to V \otimes W \otimes U,\; x\otimes y\otimes z \mapsto z \otimes y \otimes x$, using exchange
matrices in two different ways according to the following hexagon
$$
\xymatrix{
& U \otimes V \otimes W \ar[r] & V \otimes U \otimes W \ar[dr] & &\\
U \otimes W \otimes V \ar[ur] \ar[dr] & & & V \otimes W \otimes U\\
& W \otimes U \otimes V \ar[r] & W \otimes V \otimes U \ar[ur] & &
}
$$
\qed
\paragraph{Remark.} One can also deduce this proposition from Part 3.
of Proposition 2.3. Namely, one can show that if
$J(\lambda)$ is any  element of the completion of
$U(\g)\otimes U(\g)$ which satisfies  the dynamical 2-cocycle condition
(where $J_{VW}(\la)$ denotes the projection of 
$J(\la)$ to the product $V\otimes W$ of finite dimensional modules $V,W$)
then the element $R(\la)=J(\la)^{-1}J^{21}(\la)$ satisfies the 
quantum dynamical Yang-Baxter equation. 

\paragraph{Example 1.} Let us evaluate the fusion and exchange
matrices in the simplest example. Namely, take
$\g=\mathfrak{sl}_2=\C e
\oplus \C h
\oplus \C f$ and $V=\C^2=\C v_{+} \oplus \C v_{-}$ with
$$h.v_{\pm}=\pm v_{\pm}, \qquad e.v_-=v_+,\qquad e.v_+=0,\qquad f.v_-=0, \qquad f.v_+=v_-.$$
Let us compute the fusion matrix $J_{VV}(\lambda)$. By the triangularity property of $J_{VV}(\lambda)$, we have
$$
J_{VV}(\lambda)(v_{\pm} \otimes v_{\pm})=v_\pm \otimes v_\pm,\qquad
J_{VV}(\lambda) (v_- \otimes v_+)=v_- \otimes v_+,$$
so it remains to compute $J_{VV}(\lambda)(v_-\otimes v_+)$. 
 Consider the intertwiner $\Phi^{v_{-}}_\lambda: M_\lambda \to M_{\lambda+1}\otimes V$. By definition, $\Phi^{v_-}_\lambda(x_\lambda)=x_{\lambda+1} \otimes v_- + y(\lambda)fx_{\lambda+1}\otimes v_+$. To determine the function $y(\lambda)$, we use the intertwining property:
\begin{align*}
0=\Phi^{v_-}_\lambda(ex_\lambda)=(e\otimes 1+1\otimes e)\Phi^{v_-}_\lambda(x_\lambda)&=x_{\lambda+1} \otimes v_+ + y(\lambda)efx_{\lambda+1} \otimes v_+\\
&=x_{\lambda+1} \otimes v_+ + y(\lambda)(h+fe)x_{\lambda+1} \otimes v_+\\
&=x_{\lambda+1} \otimes v_+ + (\lambda+1)y(\lambda)x_{\lambda+1} \otimes v_+
\end{align*}
Hence $y(\lambda)=-\frac{1}{\lambda+1}$. It is also obvious that $\Phi^{v_+}_{\lambda+1}(x_{\lambda+1})=x_\lambda \otimes v_+$. Thus
$$\Phi^{v_+,v_-}_{\lambda}(x_\lambda)=\Phi^{v_+}_{\lambda+1}\Phi^{v_-}_\lambda(x_\lambda)=x_\lambda \otimes (v_+\otimes v_- -\frac{1}{\lambda+1} v_- \otimes v_+) + \mathrm{\;lower\;weight\;terms}.$$
Therefore $J_{VV}(\lambda)(v_+ \otimes v_-)=v_+ \otimes v_- -\frac{1}{\lambda+1}v_- \otimes v_+$, and 
$$
J_{VV}(\lambda)=\pmatrix 1&0&0&0\\0&1&0&0\\0&-\frac{1}{\lambda+1}&1&0\\0&0&0&1 \endpmatrix
$$
The exchange matrix is now easily computed. In the basis $(v_+ \otimes v_+, v_+\otimes v_-,\\v_- \otimes v_+,v_- \otimes v_-)$ it is given by
$$R_{VV}(\lambda)=\pmatrix 1 &0 & 0 & 0\\ 0&1&-\frac{1}{\lambda+1}&0\\0&\frac{1}{\lambda+1}&1-\frac{1}{(\lambda+1)^2}&0\\0&0&0&1 \endpmatrix.$$
\paragraph{2.2. Generalization to quantum groups.} The construction of intertwining operators, fusion and
exchange matrices admit natural quantum analogues. Let $U_q(\g)$ be the quantum universal enveloping algebra
associated to $\g$, as defined in \cite{CP}, Chapter 6., and for
each
$\lambda
\in
\h^*$, let
$M_{\lambda}$ be the Verma module of highest weight
$\lambda$. Then Proposition 2.2 and the definition of the fusion matrices $J_{W,V}(\lambda)$ are identical to
the classical case. In this situation, Proposition 2.3, parts 2., 3. hold. However, the fusion matrices are no longer
rational functions of $\lambda$, but rather trigonometric functions (i.e rational functions of
$q^{<\lambda,\alpha>}$, $\alpha \in \Delta$). 

Let $\mathcal{R} \in U_q(\g) \hat{\otimes} U_q(\g)$ be the universal R-matrix of $U_q(\g)$. Let $V,W$ be two finite-dimensional $U_q(\g)$-modules. The exchange matrix is defined as
$$R_{VW}(\lambda)=J_{VW}^{-1}(\lambda)\mathcal{R}^{21}_{VW}J^{21}_{W,V}(\lambda)$$
where $\mathcal{R}^{21}_{VW}$ is the evaluation of $\mathcal{R}^{21}$ on $V \otimes W$. 

In terms of intertwining operators, the exchange matrix has the following interpretation. Recall that if $V$ and $W$ are any two $U_q(\g)$-modules then $P\mathcal{R}_{VW}:V \otimes W \to W \otimes V$ is a $U_q(\g)$-intertwiner. Then $R_{VW}(\lambda) (v \otimes w)=\sum_i v_i \otimes w_i$ where $P\mathcal{R}_{WV}\Phi^{w,v}_\lambda= \sum_i \Phi_\lambda ^{v_i,w_i}$. 

With this definition, Proposition 2.4 is satisfied. The quantum analogues of the fusion and exchange matrices in
example 1 are

$$J_{VV}(\lambda)=\pmatrix
1&0&0&0\\0&1&0&0\\0&\frac{q^{-1}-q}{q^{2(\lambda+1)}-1}&1&0\\0&0&0&1
\endpmatrix,$$
$$ R_{VV}(\lambda)= \pmatrix q &0 & 0 & 0\\ 0&1&\frac{q^{-1}-q}{q^{2(\lambda+1)}-1}&0\\0&\frac{q^{-1}-q}{q^{-2(\lambda+1)}-1}&\frac{(q^{2(\lambda+1)}-q^2)(q^{2(\lambda+1)}-q^{-2})}{(q^{2(\lambda+1)}-1)^2}&0\\0&0&0&q \endpmatrix.$$

\section{The dynamical Yang-Baxter equations.}
\paragraph{3.1.} Proposition 2.4 motivates the following definition. Let $\h$ be a finite-dimensional abelian Lie
algebra and let $V$ be a semisimple $\h$-module. Let us denote by $M$ the field of meromorphic functions on
$\h^*$. Let us equip $M$ with the trivial $\h$-module structure.
\paragraph{Definition.} Let $R: V \otimes V \otimes M \to V \otimes V \otimes M$ be an $\h$-invariant and $M$-linear map. Then the \textit{quantum dynamical Yang-Baxter equation} 
(QDYBE) is the following equation with respect to $R$:
$$R^{12}(\lambda-h^{(3)})R^{13}(\lambda)R^{23}(\lambda-h^{(1)})=R^{23}(\lambda)R^{13}(\lambda-h^{(2)})R^{12}(\lambda).$$
A \textit{quantum dynamical R-matrix} is an invertible solution of this equation.
\paragraph{}It follows from Proposition 2.4 that for any simple complex Lie algebra $\g$ and for any
finite-dimensional $\g$-module $V$, the exchange matrix $R_{VV}(\lambda)$ is a quantum dynamical R-matrix.
The same is true if we replace the Lie algebra $\g$ by the quantum group $U_q(\g)$.

\paragraph{Remarks.}   
1. The usual quantum Yang-Baxter equation is recovered from the
quantum dynamical Yang-Baxter equation when $\h=0$.\\ 
2. A constant solution of the quantum dynamical Yang-Baxter
equation is the same thing as a solution of the ordinary quantum
Yang-Baxter equation which is $\h$-invariant.\\ 3. In physical
literature, the variable $\lambda$ is called a dynamical
variable. This gave rise to the name "dynamical R-matrix".

\paragraph{}Replacing $\lambda$ by $\frac{\lambda}{\gamma}$ in the
QDYBE yields the following equation
\begin{equation}\label{E:2}
\tilde{R}^{12}(\lambda-\gamma h^{(3)})\tilde{R}^{13}(\lambda) \tilde{R}^{23}(\lambda-\gamma h^{(1)})=\tilde{R}^{23}(\lambda)\tilde{R}^{13}(\lambda-\gamma h^{(2)})\tilde{R}^{12}(\lambda),
\end{equation}
which is called the quantum dynamical Yang-Baxter equation with step $\gamma$.

\begin{prop} Let $\h$ be an abelian Lie algebra. Let $V$ be a finite-dimensional semisimple $\h$-module and let ${R}:\h^* \to \mathrm{End}_\h(V \otimes V)[[\gamma]]$ be a series of meromorphic functions of the form ${R}=1-\gamma r + O(\gamma^2)$. If ${R}$ satisfies the quantum dynamical Yang-Baxter equation with step $\gamma$ then $r$ satisfies the following classical analogue of the quantum dynamical Yang-Baxter equation:
\begin{equation}\label{E:3}
\begin{split}
\sum_i &\left(x_i^{(1)} \frac{\partial r^{23}(\lambda)}{\partial x^i}-x_i^{(2)} \frac{\partial r^{13}(\lambda)}{\partial x^i} + x_i^{(3)} \frac{\partial r^{12}(\lambda)}{\partial x^i}\right) +\\
&[r^{12}(\lambda),r^{13}(\lambda)]+[r^{12}(\lambda),r^{23}(\lambda)]+[r^{13}(\lambda),r^{23}(\lambda)]=0
\end{split}
\end{equation}
where $(x_i)$ is a basis of $\h$ and $(x^i)$ is the dual basis of $\h^*$.
\end{prop}

\paragraph{}This leads to the following definition:
\paragraph{Definition.} Let $\g$ be a finite-dimensional Lie
algebra and let $\h \subset \g$ be a Lie subalgebra. The
\textit{classical dynamical Yang-Baxter equation} (CDYBE) is
equation (\ref{E:3}) with respect to a holomorphic,
$\h$-invariant function $r: \;U \to \g \otimes \g$, where $U \subset \h^*$ is an open region.  A solution to this
equation is called a \textit{classical dynamical r-matrix}.

\paragraph{Remarks.} 1. The ordinary classical Yang-Baxter equation is recovered from the classical dynamical Yang-Baxter equation when $\h=0$.\\
2. A constant solution of the classical dynamical Yang-Baxter equation is 
the same thing as an $\h$-invariant solution of the ordinary classical Yang-Baxter equation.

\paragraph{}We will now consider asymptotic behavior of
fusion and exchange matrices, and obtain solutions to the
CDYBE.  Let $\g$ be a simple
complex Lie algebra. Let $V,W$ be two finite-dimensional
$\g$-modules and let $J_{VW}(\lambda)$ and
$R_{VW}(\lambda)$ be the fusion and exchange matrices of $V$
and $W$.
\begin{prop}[\cite{EV3}] 1. The function $J_{VW}(\frac{\lambda}{\gamma})$ is regular at $\gamma=0$ for generic values of $\lambda$. 

2. Set $J_{VW}(\frac{\lambda}{\gamma})=1+\gamma j_{VW}(\lambda) + O(\gamma^2)$. Then $j_{VW}(\lambda)$ is the evaluation on $V \otimes W$ of the element 
$$j(\lambda)=-\sum_{\alpha >0}\frac{e_{-\alpha} \otimes e_\alpha}{\langle \alpha,\lambda \rangle} \in \n_- \otimes \n_+.$$
\end{prop}
\begin{cor} We have $R_{VW}(\frac{\lambda}{\gamma})=1 -\gamma r(\lambda)_{|V \otimes W}+ O(\gamma^2)$ where
\begin{equation}\label{E:2.5}
r(\lambda)=j(\lambda)-j^{21}(\lambda)=\sum_{\alpha>0}\frac{e_\alpha \otimes e_{-\alpha} -e_{-\alpha} \otimes e_\alpha}{\langle \alpha, \lambda \rangle}.
\end{equation}
\end{cor}

A proof of Proposition 3.2, which is based on computing the asymptotics 
of intertwining operators at $\lambda\to\infty$, is given in \cite{EV3}. 
Later we will give another proof of this Proposition.
 
\paragraph{}It follows from Proposition 3.1 that $r(\lambda)$  in (\ref{E:2.5}) is a classical dynamical r-matrix.
Let us call it the \textit{basic rational dynamical r-matrix}.

\paragraph{}Proposition 3.2 and Corollary 3.1 have natural quantum analogues. Let $U_q(\g)$ be the quantum
group associated to $\g$ with quantum parameter $q=e^{-\varepsilon \gamma /2}$ for some fixed $\varepsilon
\in \C$ and formal parameter $\gamma$. Let $V,W$ be two finite-dimensional $U_q(\g)$-modules and let
$R_{VW}(\lambda)$ be the exchange matrix. Set $\tilde{R}_{VW}(\lambda)=R_{VW}(\frac{\lambda}{\gamma})$.
\begin{prop}[\cite{EV3}] We have
$$\tilde{R}_{VW}(\lambda)=1-\gamma r^\varepsilon_{VW}(\lambda) + O(\gamma^2)$$
where $r^\varepsilon_{VW}: \h^* \to \mathrm{End}_\h(V \otimes W)$ is the evaluation on $V \otimes W$ of the following universal element:
\begin{equation}\label{E:4}
r^\varepsilon(\lambda)=\frac{\varepsilon}{2}\Omega + \sum_{\alpha >0} \frac{\varepsilon}{2} \mathrm{cotanh} \left( \frac{\varepsilon}{2} \langle \alpha,\lambda\rangle \right) (e_\alpha \otimes
e_{-\alpha} - e_{-\alpha} \otimes e_\alpha) \in \g \otimes \g
\end{equation}
where $\Omega \in S^2\g$ is the inverse element to the form $(\;,\;)$ (the \textit{Casimir} element).
\end{prop}
It follows from Proposition 3.1 that $r^\varepsilon(\lambda)$ is a
solution of the CDYBE. Let us
call it the \textit{basic trigonometric dynamical r-matrix}.
\paragraph{3.2. Quantization and quasiclassical limit.} Let $\h$
be an abelian Lie algebra and let $V$ be a finite-dimensional
semisimple $\h$-module. Let $r: \h^* \to \mathrm{End}_\h(V
\otimes V)$ be a classical dynamical r-matrix. Suppose that
$R:\h^* \to \mathrm{End}_\h(V \otimes V)[[\gamma]]$ is of the
form $R=1-\gamma r + O(\gamma^2)$ and satisfies the QDYBE.
\paragraph{Definition.} $R$ is called a \textit{quantization} of $r$. Conversely, $r$ is called the \textit{quasiclassical limit} of $R$.
\paragraph{}For instance, the exchange matrix $\tilde{R}_{VV}(\lambda)$ constructed from a Lie algebra $\g$ is a quantization of the evaluation on $V \otimes V$ of the basic rational dynamical r-matrix. Similarly, exchange matrices constructed from  quantum groups provide quantization of the basic trigonometric dynamical r-matrix.
\paragraph{Conjecture:} Any classical dynamical $r$-matrix admits a quantization.
\paragraph{}Notice that when $\h=0$, this conjecture reduces to the conjecture of Drinfeld \cite{Dr} about quantization of classical (non-dynamical) r-matrices, which was proved in \cite{EK}. In the skew-symmetric case, 
the conjecture was recently proved in \cite{Xu2,Xu3}
(under some minor technical assumptions), using the theory of Fedosov
quantization. Also, the conjecture is proved in \cite{ESSch} 
for the classical dynamical r-matrices on simple Lie algebras
classified in \cite{S}. 

\paragraph{3.3. Unitarity conditions.} Recall the following notions introduced by Drinfeld. A classical r-matrix $r \in \g \otimes \g$ is a \textit{quasitriangular} structure on a Lie algebra $\g$ if $r+r^{21} \in (S^2\g)^\g$. It is a \textit{triangular} structure on $\g$ if $r+r^{21}=0$. 

This definition is natural in the theory of Lie bialgebras. Namely, a classical r-matrix $r \in \g \otimes \g$
defines a Lie bialgebra structure on $\g$ by $\delta: \g \to \Lambda^2 \g,\; x \mapsto [1 \otimes x + x \otimes
1,r]$ if and only if $r+r^{21} \in (S^2\g)^\g$. In the case of a simple Lie algebra $\g$ we have
$(S^2\g)^\g=\C\Omega$, so that a classical r-matrix $r$ is quasitriangular if $r+r^{21}=\varepsilon \Omega$
for some $\varepsilon \in \C$, and it is triangular if moreover
$\varepsilon=0$. 
\paragraph{}This leads one to make the following definition:
\paragraph{Definition.} A classical dynamical r-matrix $r: \h^* \to (\g \otimes \g)^\h$ has \textit{coupling
constant}
$\varepsilon$ if
\begin{equation}\label{E:5}
r+r^{21}=\varepsilon \Omega.
\end{equation}
Equation (\ref{E:5}) is called the \textit{unitarity condition}. Notice that the basic rational dynamical r-matrix
$r(\lambda)$ and the basic trigonometric dynamical r-matrix $r^\varepsilon (\lambda)$ have coupling constants
$0$ and $\varepsilon$ respectively.

\section{Classification of classical dynamical $r$-matrices.}
\paragraph{}In this section, we give the classification of all dynamical r-matrices $r:\h^* \to \g \otimes \g$ which have coupling constant $\varepsilon \in \C$.
\paragraph{4.1. Gauge transformations.} Consider the following operations on meromorphic maps $r: \h^* \to (\g \otimes \g)^\h$.
\begin{enumerate}
\item $r(\lambda) \mapsto r(\lambda) + \sum_{i<j} C_{ij}(\lambda) x_i \wedge x_j,$,
where $\sum_{i,j} C_{ij}(lambda) d\lambda_i \wedge
d\lambda_j$ is a closed meromorphic 2-form.
\item $r(\lambda) \mapsto r(\lambda-\nu),$  
where $\nu \in \h^*$.
\item
$r(\lambda) \mapsto (A \otimes A) r(A^*\lambda),$ 
where $A\in W$, the Weyl group of $\g$.
\end{enumerate}
\begin{lem} Transformations 1-3 preserve the set of classical dynamical r-matrices. \end{lem}
\noindent The proof is straightforward. 

Two classical dynamical r-matrices which can be obtained one from the other by a sequence of such transformations will be called \textit{gauge-equivalent}.
\paragraph{4.2. Classification of dynamical r-matrices with zero coupling constant.}  Let  $\l \supset \h$ be a reductive Lie subalgebra of $\g$. Define
\begin{equation}\label{E:2.65}
r^\l(\lambda)=\sum_{\underset{e_\alpha \in \l}{\alpha>0}} \frac{e_\alpha \otimes e_{-\alpha} - e_{-\alpha} \otimes e_{\alpha}}{(\lambda,\alpha)}.
\end{equation}
It is clear that this is the image of the basic rational dynamical r-matrix of $\l$ under the embedding $\l \subset \g$.
\begin{theo}[\cite{EV1}] Any classical dynamical r-matrix $r: \h^* \to (\g \otimes \g)^\h$ with zero coupling constant is gauge-equivalent to $r^\l(\lambda)$ for some $\l$.
\end{theo}
\paragraph{4.3. Classification of dynamical r-matrices with coupling constant $\varepsilon \in \C^*$.} Let $X \subset \Pi$, and denote by $\langle X \rangle \subset \Delta$ the set of all roots which are linear combinations of elements in $X \cup -X$. For any $\alpha \in \Delta$ introduce a meromorphic function $\varphi_\alpha: \h^* \to \C$ by the following rule. Set $\varphi_\alpha(\lambda)=\frac{\varepsilon}{2}$ if $\alpha \in \Delta^+\backslash\langle X \rangle$, $\varphi_\alpha(\lambda)=-\frac{\varepsilon}{2}$ if $\alpha \in \Delta^-\backslash\langle X \rangle$ and 
$$
\varphi_\alpha (\lambda)= \frac{\varepsilon}{2} \mathrm{cotanh} \left( \frac{\varepsilon}{2}(\lambda,\alpha)\right)
$$ 
if $\alpha \in \langle X \rangle$.
\begin{theo}[\cite{EV1}] Let $X \subset \Pi$. Set  
$$r_X(\lambda)=\frac{\varepsilon}{2}\Omega + \sum_{\alpha \in \Delta} \varphi_\alpha (\lambda) e_\alpha \otimes e_{-\alpha}.$$
Then $r^\varepsilon_X(\lambda)$ is a classical dynamical r-matrix with coupling constant $\varepsilon$. Moreover, any classical dynamical r-matrix with coupling constant $\varepsilon$ is gauge-equivalent to $r^\varepsilon_X(\lambda)$ for a suitable $X \subset \Pi$.
\end{theo}
\paragraph{Remarks.} 1. The basic trigonometric dynamical r-matrix $r^\varepsilon(\lambda)$ is obtained when we take $X=\Pi$. Moreover, the r-matrix $r^\varepsilon_X(\lambda)$ is equal to a
limit of $r^\varepsilon(\lambda-\nu)$ when $\nu$ tends to infinity in $\h^*$ in an appropriate direction.
In other words, every classical dynamical r-matrix with nonzero coupling constant $\varepsilon$ is a limiting
case of the basic trigonometric r-matrix.\\ 2. Let $W$ be the Weyl group of $\g$, and let $w \in W$. Let
$\lambda \in
\h^*$ tend to infinity in a generic way in the Weyl chamber associated to $w$. Then
$$\mathrm{lim}\; r^1(\lambda)=\frac{1}{2}\sum_{i} x_i \otimes x_i + \sum_{\alpha \in w(\Delta^+)} e_\alpha \otimes e_{-\alpha},$$
which is the standard classical r-matrix corresponding to the polarization of $\g$ associated to $w$. Hence the basic trigonometric dynamical r-matrix $r^1(\lambda)$ interpolates all $\h$-invariant classical (non-dynamical) r-matrices $r$ satisfying $r+r^{21}=\Omega$, (up to the addition of a skew 2-form in $\Lambda^2 \h$).
\paragraph{}A classification of all classical dynamical r-matrices $r:\l^* \to (\g \otimes \g)^\l$ where $\g$
is a simple Lie algebra and $\l\subset \h$ is given in \cite{S}. This classification generalizes both the above
classification (when $\l=\h$) and the Belavin-Drinfeld classification of classical r-matrices (when $\l=0$) (see
Appendix A). 

\section{Classical dynamical r-matrices and Poisson-Lie groupoids}
\paragraph{} In this section we give a geometric interpretation of the
CDYBE. Let us first briefly recall the relationship between
the classical Yang-Baxter equation and the theory of
Poisson-Lie groups, developed by Drinfeld.
\paragraph{5.1. Poisson-Lie groups.} Let $G$ be a (complex or real) Lie group, let $\g$ be its Lie algebra  and let
$\mathcal{O}(G)$ be the algebra of regular functions on $G$. Let $\{\,,\,\}: \mathcal{O}(G) \times \mathcal{O}(G)
\to \mathcal{O}(G)$ be a Poisson structure on $G$. Let $\Pi$ be the Poisson bivector field, defined by the
relation $\{f,g\}=df \otimes dg (\Pi)$. Recall that $(G,\{\,,\,\})$ is called a Poisson-Lie group if the
multiplication map
$m: G \times G \to G$ is a Poisson map.\\
\hbox to1em{\hfill}Let $\rho \in \Lambda^2 \g$ and consider the following bivector field:
$$\Pi_\rho=R_\rho-L_\rho,$$
where $R_\rho$ (resp. $L_\rho$) is the left-invariant (resp. right-invariant) bivector field satisfying $(R_\rho)_{|e}=\rho$ (resp. $(L_\rho)_{|e}=\rho$);
in other words, $R_\rho,L_\rho$ stand for the translates of $\rho$ by right and left shifts respectively.
\begin{prop}[Drinfeld] The bivector $\Pi_\rho$ defines a Poisson-Lie group structure on $G$ if and only if
$$[\rho_{12},\rho_{13}]+[\rho_{13},\rho_{23}]+[\rho_{12},\rho_{23}]\in (\Lambda^3 \g)^\g.$$
\end{prop}
When this is the case, $G$ is called a {\textit{coboundary}} Poisson-Lie group. Two cases are of special interest:
\begin{enumerate}
\item The exists $T \in (S^2\g)^\g$ such that $[\rho_{12},\rho_{13}]+[\rho_{13},\rho_{23}]+[\rho_{12},\rho_{23}]=\frac{1}{4}[T_{12},T_{23}]$. This implies that $r=\rho + \frac{1}{2}T$ satisfies the classical Yang-Baxter equation, and $\Pi_\rho=R_r-L_r$. In this case $G$ is called a {\textit{quasitriangular}} Poisson-Lie group.
\item We have $[\rho_{12},\rho_{13}]+[\rho_{13},\rho_{23}]+[\rho_{12},\rho_{23}]=0$. In this case, $G$ is called a {\textit{triangular}} Poisson-Lie group.
\end{enumerate}
\paragraph{5.2. Poisson-Lie groupoids.} It turns out that, in order to generalize this correspondence to the dynamical case, groups must be replaced by groupoids. Recall that a groupoid is a (small) category where all morphisms are isomorphisms. It is equivalent to the following data: two sets $X$ and $P$ (the set of morphisms, or the groupoid itself,  and the set of objects, or the base, respectively), two surjective maps $s,t: X \to P$ (the source and target maps), an injective map $E:P \to X$ (the identity morphisms), a map
$m:\{(a,b) \in X \times X\;|\; t(a)=s(b)\} \to X$ ($m(a,b)=b \circ a$, the composition of morphisms), and an involution $i:X \to X$ such that $s(i(x))=t(x),\;t(i(x))=s(x),\;m(i(x),x)=\mathrm{Id}_{s(x)}$ and $m(x,i(x))=\mathrm{Id}_{t(x)}$ for all $x \in X$, satisfying some obvious axioms. One can visualize elements of $X$ as arrows $s(a) \stackrel{a}{\to} t(a)$. \\
 Note that when $|P|=1$, the notion of a groupoid coincides with the notion of a group.\\
\hbox to1em{\hfill} A Lie groupoid is a groupoid with a smooth structure (in particular, the sets of objects and the sets of morphisms are smooth manifolds and the structure maps are smooth, see \cite{M}).
\paragraph{}Now we would like to generalize the notion of a
Poisson-Lie group to groupoids. The usual definition does not
generalize directly since if $X$ is a Lie groupoid and a Poisson manifold 
then the set of points $(a,b) \in X^2$ for
which the multiplication is defined is not necessarily a Poisson submanifold,
so we cannot require that the multiplication map be Poisson. But
this difficulty can be bypassed using the following
observation:
\begin{prop} Let $X,Y$ be two Poisson manifolds and let $f: X \to Y$ be a smooth map. Consider the graph $\Gamma_f=\{(x,f(x))\} \subset X \times \overline{Y}$, where $\overline{Y}$ is the manifold $Y$, with the opposite Poisson structure $\{\,,\,\}_{\overline{Y}}=-\{\,,\,\}_Y$. Then $f$ is a Poisson map if and only if $\Gamma_f$ is a coisotropic submanifold of $X \times \overline{Y}$, i.e if and only if for any $z \in \Gamma_f$, $(T_z\Gamma_f)^\perp \subset T_z^*( X \times \overline{Y})$ is an isotropic subspace with respect to the Poisson
form 
$\Pi$ on $T_z(X\times\overline{Y})^*$.
\end{prop}
\paragraph{}This gives rise to the following notion of a Poisson-Lie groupoid, first introduced by Weinstein \cite{W}.
\paragraph{Definition.} A Lie groupoid $X$ with a Poisson structure is called a \textit{Poisson-Lie groupoid} if $\Gamma_m \subset X \times X \times \overline{X}$ is a coisotropic submanifold.
\paragraph{}We now restrict ourselves to a particular class of
Lie groupoids. Let $G$ be a Lie group, let $\g$ be its Lie algebra,
$\h \subset \g$ a subalgebra and $H$ a Lie subgroup of $G$ with
Lie algebra $\h$. Let $U
\subset
\h^*$ be an open set.  Consider the following groupoid: $X=U \times G
\times U$, $P=U$ with $s(u_1,g,u_2)=u_1$, $t(u_1,g,u_2)=u_2$. The
composition $m((u_1,g,u_2),(u_3,g',u_4))$ is defined only when
$u_2=u_3$ and $m((u_1,g,v),(v,g',u_4))=(u_1,gg',u_4)$. If $a$ is
a function on $U$ we set $a_1=s^*(a) \in \mathcal{O}(X)$ and
$a_2=t^*(a) \in \mathcal{O}(X)$. Let $\rho: U \to \Lambda^2\g$
be a regular function. \\
\hbox to1em{\hfill}The group $H^2$ acts on $X$ by
$$(h_1,h_2)(u_1,g,u_2)=(\mathrm{Ad}^*(h_1)u_1,h_1gh_2^{-1},\mathrm{Ad}^*(h_2)u_2).$$
We want to define a Poisson structure on $X$ for which
$(-s,t)$ is a moment map for this action. This forces the
following relations
\begin{equation}\label{E:21}
\begin{split}
\{a_1,b_1\}=-[a,b]_1,\qquad &\{a_2,b_2\}=[a,b]_2,\qquad \{a_1,b_2\}=0,\\
\{a_1,f\}=&R_a f,\qquad \{a_2,f\}=L_af.
\end{split}
\end{equation}
We try to complete the definition of the Poisson structure on $X$ by adding
the relation
\begin{equation}\label{E:22}
\{f,g\}=(df \otimes dg)(R_{\rho(u_1)}-L_{\rho(u_2)}),
\end{equation}
where $f,g$ are any functions on $X$ pulled back from $G$ and $a,b$ are linear functions on $U$.
\begin{prop}[\cite{EV1}] Formulae (\ref{E:21}) and (\ref{E:22}) define a Poisson-Lie groupoid structure on $X$ if and only if
\begin{enumerate}
\item
$$
\sum_i \left(x_i^{(1)} \frac{\partial \rho^{23}}{\partial x^i}-x_i^{(2)} \frac{\partial \rho^{13}}{\partial x^i} + x_i^{(3)} \frac{\partial \rho^{12}}{\partial x^i}\right) +
[\rho^{12}, \rho^{13}]+[\rho^{12},\rho^{23}]+[\rho^{13},\rho^{23}]$$
is a constant $\g$-invariant element of $\Lambda^3\g$, and 
\item $\rho$ is $\h$-invariant.
\end{enumerate}
\end{prop}

When this is the case, $X$ is called a {\textit{coboundary}} dynamical Poisson-Lie groupoid, which will be denoted by $X_r$. 
Two cases are of special interest:
\begin{enumerate}
\item The exists $T \in (S^2\g)^\g$ such that 
\begin{equation*}
\begin{split}
\frac{1}{4}[T_{12},T_{23}]=\sum_i &\left(x_i^{(1)} \frac{\partial
r^{23}(\lambda)}{\partial x^i}-x_i^{(2)} \frac{\partial
r^{13}(\lambda)}{\partial x^i} + x_i^{(3)} \frac{\partial
r^{12}(\lambda)}{\partial x^i}\right)\\
 &+ [\rho_{12},\rho_{13}]+[\rho_{13},\rho_{23}]+[\rho_{12},\rho_{23}].
\end{split}
\end{equation*}
 This implies that $r=\rho + \frac{1}{2}T$ satisfies the classical dynamical Yang-Baxter equation. In this case
$X_r$ is called a {\textit{quasitriangular}} dynamical Poisson-Lie
groupoid.
\item We have
\begin{equation*}
\begin{split}
\sum_i \bigg(x_i^{(1)} \frac{\partial
r^{23}(\lambda)}{\partial x^i}&-x_i^{(2)} \frac{\partial
r^{13}(\lambda)}{\partial x^i} + x_i^{(3)} \frac{\partial
r^{12}(\lambda)}{\partial x^i}\bigg)\\
 &+ [\rho_{12},\rho_{13}]+[\rho_{13},\rho_{23}]+[\rho_{12},\rho_{23}]=0.
\end{split}
\end{equation*}
In this case, $X_r$ is called a
{\textit{triangular}} dynamical Poisson-Lie groupoid.
\end{enumerate}

Thus, the basic rational solution defined above gives rise 
to a triangular dynamical Poisson-Lie groupoid, and the basic trigonometric 
solution gives rise to a quasitriangular one. 
        
\section{Classification of quantum dynamical R-matrices}

In this section we give the classification of all quantum dynamical $R$-matrices $R: \h^* \to \mathrm{End}_\h(V \otimes V)$, where $\h$ is the Cartan subalgebra of $\mathfrak{gl}(n,\C)$ consisting of diagonal matrices, and $V=\C^n$ is the vector representation, which satisfy an additional
\textit{Hecke condition}, a quantum analogue of the unitarity condition.   
\paragraph{6.1. Hecke condition.} Let $\h$ be the abelian Lie algebra 
of diagonal N by N matrices,
and let $V$ be the standard N-dimensional $\h$-module.
Let $h_1,\ldots h_n$ be the standard basis of $\h$, $\lambda_1,\ldots \lambda_n$ be the corresponding coordinate functions on $\h^*$, and $V_i,\; i=1,\ldots n$ be the (one-dimensional) weight subspaces of $V$ of weight $\omega_i$ where $\langle \omega_i,h_j \rangle=\delta_{ij}$. 

Consider the $\h$-module $V\otimes V$. Its weight subspaces are 
$V_a\otimes V_b\oplus V_b\otimes V_a$ and $V_a \otimes V_a$. 

\paragraph{Definition.} An operator $R: \h^* \to \mathrm{End}_\h(V \otimes V)$ satisfies the \textit{Hecke condition} with parameter $q \in \C^*$ if
the eigenvalues of $PR$ (where $P$ is the permutation matrix) are 
$1$ on $V_a\otimes V_a$ and $1,-q$ on 
$V_a\otimes V_b\oplus V_b\otimes V_a$. 

\paragraph{}This condition can be thought of as a quantum analogue of the unitarity condition for classical r-matrices, since it is easy to show that 
the quasiclassical limit of an operator satisfying the Hecke condition
satisfies the unitarity condition.
 In particular if $R$ satisfies the Hecke condition with $q=1$ then $RR^{21}=1$, which can be thought of as a quantization of the relation $r+r^{21}=0$.\\
\hbox to1em{\hfill}The terminology comes from the following remark: if $R$ is a $\lambda$-independent solution of the quantum dynamical Yang-Baxter equation satisfying the Hecke condition with parameter $q$ then $\check{R}$ defines a representation of the Hecke algebra $H_p$ of type $A_{p-1}$ on the space $V^{\otimes p}$ for any $p>1$. Similar representations can be defined 
for dynamical R-matrices (see \cite{EV2} and Section 7).

\paragraph{6.2. Gauge transformations.} 
Let $R(\lambda)$ be a quantum dynamical R-matrix satisfying Hecke 
condition with parameter $q$. The weight-zero and Hecke conditions imply 
that
\begin{equation}\label{E:formulaforR}
R(\lambda)=\sum_{a} E_{aa}
\otimes E_{aa} +
\sum_{a\neq b}
\alpha_{ab}(\lambda) E_{aa} \otimes E_{bb} + \sum_{a \neq b}
\beta_{ab}(\lambda) E_{ab} \otimes E_{ba}
\end{equation}
 where $E_{ij}$ is the
elementary matrix,  and $\alpha_{ab},\beta_{ab}$ are meromorphic
functions $\h^* \to \C$. So it is enough to look for solutions of this
form. 

As in the classical case, we will give the classification 
of solutions up to some group of transformations.
\paragraph{Definition.} A \textit{multiplicative 2-form} on $V$ is a collection meromorphic functions $\{\varphi_{ab}: \h^* \to \C\}_{a,b=1}^n$ satisfying $\varphi_{ab}\varphi_{ba}=1$ for all $a,b$. A multiplicative 2-form $\{\varphi_{ab}(\lambda)\}$ is \textit{closed} if for all $a,b,c$,
$$\frac{\varphi_{ab}(\lambda)}{\varphi_{ab}(\lambda-\omega_c)}\frac{\varphi_{bc}(\lambda)}{\varphi_{bc}(\lambda-\omega_a)}\frac{\varphi_{ca}(\lambda)}{\varphi_{ca}(\lambda-\omega_b)}=1.$$
\paragraph{}Consider the following operations on meromorphic weight-zero maps $R: \h^* \to \mathrm{End}_\h(V \otimes V)$
of the form (\ref{E:formulaforR}):
\begin{enumerate}
\item $$R(\lambda) \mapsto \sum_{a}  E_{aa}
\otimes E_{aa} + \sum_{a\neq b}
\varphi_{ab}(\lambda)\alpha_{ab}(\lambda) E_{aa} \otimes E_{bb} +
\sum_{a \neq b} \beta_{ab}(\lambda) E_{ab} \otimes E_{ba},$$ where
$\{\varphi_{ab}(\lambda)\}$ is a closed multiplicative 2-form on
$V$,
\item $R(\lambda) \mapsto R(\lambda-\nu)$ where $\nu$ is a pseudoconstant, 
i.e. a meromorphic function $\h^*\to \h^*$ such that 
$\nu(\lambda+\omega_i)=\nu(\lambda)$ for all $i$ (for example, a constant),

\item $R(\lambda) \mapsto (\sigma\otimes \sigma) 
R(\sigma^{-1}\lambda)(\sigma^{-1}\otimes \sigma^{-1})$ 
where $\sigma \in \mathfrak{S}_n$ acts on $V$ and $\h^*$ by permutation of coordinates.
\end{enumerate}
\paragraph{Remark.} Here we allow to perform transformation 2 only if 
the answer is meromorphic. 
\begin{lem} Tranformations $1-3$ preserve the set of quantum dynamical R-matrices.
\end{lem}

Two R-matrices which can be obtained one from the other by a sequence of such transformations are said to be \textit{gauge-equivalent}.
\paragraph{6.3. Classification for $q=1$.} Let $X$ be a subset of $\{1,\ldots n\}$ and write $X=X_1 \cup \ldots \cup X_k$ where $X_i=\{a_i\ldots b_i\}$ are disjoint intervals. Set
$$R_X(\lambda)=\sum_{a,b=1}^n E_{aa} \otimes E_{bb} + \sum_{l=1}^k \sum_{\underset{a \neq b}{a,b \in X_l}}\frac{1}{\lambda_a-\lambda_b}(E_{aa} \otimes E_{bb} + E_{ba}\otimes E_{ab}).$$
\begin{theo}[\cite{EV2}] Let $X\subset \{1,\ldots n\}$. 
Then $R_X(\lambda)$ is a quantum dynamical R-matrix satisfying the Hecke condition with $q=1$. Moreover, any dynamical R-matrix $R:\h^* \to \mathrm{End}_\h(V \otimes V)$ is gauge-equivalent to $R_X(\lambda)$ for a unique subset $X \subset \{1,\ldots n\}$.\end{theo}
\paragraph{Remark.} The function $R_X(\lambda/\gamma)$ is, up to a gauge transformation,  a quantization in the sense of Section
3 of the rational classical dynamical r-matrix (\ref{E:2.65}) corresponding to the reductive subalgebra of
$\mathfrak{gl}(n)$ spanned by root subspaces $\g_\alpha, \g_{-\alpha}$ for $\alpha \in X$. 

The most interesting solution $R_X$ corresponds to the case when 
$X=\{1,...,n\}$. We will call it {\it the basic rational solution}
of the QDYBE.  

\paragraph{6.4. Classification for  $q\neq1$.} Let $\varepsilon \not\in 2i\pi\Z$ and set $q=e^\varepsilon$. Let $X$ be a subset of $\{1,\ldots n\}$ and again write $X=X_1 \cup \ldots \cup X_k$ where $X_i=\{a_i\ldots b_i\}$ are disjoint intervals. Set
$$R^\varepsilon_X(\lambda)=\sum_{a} E_{aa} \otimes E_{aa} + \sum_{a\ne b}\alpha_{ab}(\lambda) E_{aa} \otimes E_{bb} + \sum_{a \neq b} \beta_{ab}(\lambda) E_{ab} \otimes E_{ba},$$
where $\alpha_{ab}(\lambda)=q+\beta_{ab}(\lambda)$ and where $\beta_{ab}(\lambda)$ is defined as follows: $\beta_{ab}=\frac{q-1}{q^{\lambda_a-\lambda_b}-1}$ if $a,b \in X_l$ for some $1 \leq l\leq k$, $\beta_{ab}(\lambda)=1-q$ otherwise if $a>b$ and $\beta_{ab}(\lambda)=0$ otherwise if $a<b$.

\begin{theo}[\cite{EV2}] Let $X\subset \{1,\ldots n\}$. Then $R^\varepsilon_X(\lambda)$ is a quantum dynamical R-matrix satisfying the Hecke condition with $q=e^\varepsilon$. Moreover, any dynamical R-matrix $R:\h^* \to \mathrm{End}_\h(V \otimes V)$ is gauge-equivalent to $R_X(\lambda)$ for a unique subset $X \subset \{1,\ldots n\}$.\end{theo}
\paragraph{Remark.} It can be checked that the $R^\varepsilon_X(\lambda)$ yield (again up to gauge
transformations) quantizations of the trigonometric classical dynamical r-matrices with coupling constant
$\varepsilon$ appearing in Theorem 4.2.

The most interesting solution $R^\varepsilon_X$ corresponds to the case when 
$X=\{1,...,n\}$. We will call it {\it the basic trigonometric solution}
of the QDYBE.  

\paragraph{6.5. The fusion and exchange matrices for the vector representation of 
classical and quantum $gl_n$.}
The above classification can be applied to compute 
the fusion and exchange matrices for the vector representation. 
Namely, we have:
\begin{theo}[\cite{EV3}] 1. Let $\g=\mathfrak{gl}_n$ and let $V=\C^n$ be the vector representation. Then
\begin{align*}
J_{VV}(\lambda)&= 1+\sum_{a<b} \frac{1}{\lambda_b-\lambda_a+a-b}E_{ba} \otimes E_{ab}\\
R_{VV}(\lambda)&=\sum_{a=1}^n E_{aa}\otimes E_{aa} + \sum_{a\neq b} \frac{1}{\lambda_a-\lambda_b+b-a} E_{ba} \otimes E_{ab} + \sum_{a<b} E_{aa} \otimes E_{bb}\\
&\qquad -\sum_{a>b}\frac{(\lambda_b-\lambda_a+a-b-1)(\lambda_b-\lambda_a+a-b+1)}{(\lambda_b-\lambda_a+a-b)^2}E_{aa}\otimes E_{bb}.
\end{align*}
2. Let $V=\C^n$ be the representation of $U_q({\frak gl}_N)$ which is 
the q-analog of the vector representation. Then
\begin{align*}
J_{VV}(\lambda)&= 1+\sum_{a<b} \frac{q^{-1}-q}{q^{2(\lambda_a-\lambda_b+b-a)}-1}E_{ba} \otimes E_{ab}\\
R_{VV}(\lambda)&=q\sum_{a=1}^n E_{aa}\otimes E_{aa} + \sum_{a\neq b} \frac{q^{-1}-q}{q^{2(\lambda_a-\lambda_b+b-a)}-1} E_{ba} \otimes E_{ab} + \sum_{a<b} E_{aa} \otimes E_{bb}\\
&\qquad +\sum_{a>b}\frac{(q^{2(\lambda_b-\lambda_a+a-b)}-q^{-2})(q^{2(\lambda_b-\lambda_a+a-b)}-q^{2})}{(q^{2(\lambda_b-\lambda_a+a-b)}-1)^2} E_{aa}\otimes E_{bb}.
\end{align*}
\end{theo}
\noindent
\textit{Proof.} The proof relies on explicit computations and on the classification of quantum dynamical
R-matrices (Theorems 6.1 and 6.2). More precisely, it is possible to compute explicitly the coefficients of $J$
correpsonding to simple roots, and all the other coefficients are then uniquely determined by Theorems 6.1 and
6.2. \qed

\paragraph{Remark.} 
The matrix coefficients of $J_{VV}(\lambda)$ for nonsimple roots 
are not as easily computed directly as those for simple roots. 
The above approach allows one to avoid this calculation. 

\section{Quantum dynamical R-matrices and quantum groupoids}

In this section we will give a "noncommutative geometric" 
interpretation of the QDYBE which is analogous to
the  geometric interpretation of the 
CDYBE given above. More
precisely,  to solutions of the QDYBE
we will associate,  following \cite{F},\cite{EV2}, a kind of
quantum group, more precisely a 
\textit{Hopf algebroid} (or \textit{quantum groupoid}). 

The general notion of a Hopf algebroid was introduced in \cite{Lu}. 
However, here it will be sufficient to use a less general notion, 
that of an $H$-Hopf algebroid, which was introduced in \cite{EV2}.
Our exposition will follow \cite{EV2,EV3}.

\paragraph{7.1. $H$-bialgebroids.}Let $H$ be a commutative and cocommutative 
finitely generated Hopf algebra over $\C$,
$T=\mathrm{Spec}\; H$ the corresponding commutative affine
algebraic group. Assume that $T$ is connected. 
Let $M_{T}$ denote the field of meromorphic functions on $T$.
Let us introduce the following definitions. 
\paragraph{Definition.} {\it{An $H$-algebra }} is an associative algebra
$A$ over $\C$ with $1$, endowed with an $T$-bigrading
$A=\oplus_{\al,\beta\in T}A_{\al\beta}$
(called the weight decomposition), and two algebra embeddings
$\mu_l,\mu_r:M_{T}\to A_{00}$ (the left and the right moment maps), such
that for any $a\in A_{\al\beta}$ and $f\in M_{T}$, we have
\begin{equation}\label{mmm}
\mu_l(f(\la))a=a\mu_l(f(\la+\al)),\quad
\mu_r(f(\la))a=a\mu_r(f(\la+\beta)).
\end{equation}\\
\hbox to1em{\hfill}{\it A morphism} $\phi:A\to B$ of two $H$-algebras
is an algebra homomorphism, preserving
the moment maps. 

\paragraph{Example 1.} Let $D_T$ be the algebra of difference operators
$M_{T}\to M_{T}$, i.e.
the operators of the form $\sum_{i=1}^nf_i(\la)\Ti_{\beta_i}$, where
$f_i\in M_{T}$, and
for $\beta\in T$ we denote by $\Ti_\beta$ the field automorphism
of $M_{T}$ given by $(\Ti_\beta f)(\la)=f(\la+\beta)$.

The algebra $D_T$ is an  example of an $H$-algebra  
if we define the weight decomposition by $D_T=\oplus (D_T)_{\al\beta}$,
where $(D_T)_{\al\beta}=0$ if $\al\ne\beta$, and
$(D_T)_{\al\al}=\{f(\la)\Ti_\al^{-1}: f\in M_{T}\}$, and the moment maps
$\mu_l=\mu_r:M_{T}\to (D_T)_{00}$
to be the tautological isomorphism.

\paragraph{Example 2.} This is a generalization of Example 1. 
Let $W$ be a diagonalizable $H$-module, 
$W=\oplus_{\la \in T} W[\la], \, W[\la]=\{w\in W\,|\, aw=\la(a)w, \text{for all}\,
a\in H \}$,  and let $D^\al_{T, W}\subset
\text{Hom}_\C(W,W\T D_T)$ be the space of all difference
operators on $T$ with coefficients in $\End_\C(W)$, which have weight $\al \in T$ with
respect to the action of $H$ in $W$.

Consider the algebra $D_{T, W}=\oplus_\al D_{T, W}^\al$.
This algebra has a weight decomposition
$D_{T,W}=\oplus_{\al,\beta} (D_{T,W})_{\al\beta}$ defined as follows:
if $g\in \text{Hom}_\C(W,W\T M_{T})$ is
an operator of weight $\beta-\al$,
then $g\Ti_{\beta}^{-1}\in (D_{T,W})_{\al\beta}$.

Define the moment maps $\mu_l,\mu_r: M_{T}\to (D_{T, W})_{00}$
by the formulas $\mu_r(f(\la))=f(\la)$,
$\mu_l(f(\la))=f(\la- h)$ where $f(\la - h)w=f(\la - \mu)w$
if $w\in W[\mu],\, \mu \in T$.  The algebra $D_{T, W}$ equipped with this weight
decomposition and these moment maps is an $H$-algebra.

\paragraph{}Now let us define the tensor product of $H$-algebras. 
Let $A,B$ be two $H$-algebras
and $\mu_l^A,\mu_r^A,\mu_l^B,\mu_r^B$ their moment maps.
Define their {\it matrix tensor product}, 
$A\wo B$, which is also an $H$-algebra.
Let
\begin{equation}\label{pr}
(A\wo B)_{\al\delta}:=\oplus_{\beta}A_{\al\beta}
\T_{M_{T}} B_{\beta\delta},
\end{equation}
 where $\T_{M_{T}}$ means
the usual tensor product
modulo the relation $\mu_r^A(f)a\T b=a\T \mu_l^B(f)b$, for any
$a\in A,b\in B, f\in M_{T}$.
Introduce a multiplication in $A\wo B$ by the rule $(a\T b)(a'\T b')=
aa'\T bb'$. It is easy to check that the multiplication is well defined.
 Define
the moment maps for $A\wo B$ by
$\mu_l^{A\wo B}(f)=\mu_l^A(f)\T 1$,
$\mu_r^{A\wo B}(f)=1\T \mu_r^B(f)$.

For any $H$-algebra $A$, the algebras
$A\wo D_T$ and $D_T \wo A$ are canonically isomorphic to $A$.
In particular, $D_T$ is canonically isomorphic to $D_T\wo D_T$. 
Thus the category of $H$-algebras
equipped with the product $\wo$ is a
monoidal category, where the unit object is $D_T$.

\paragraph{}Now let us define the notions of a coproduct and a counit on an $H$-algebra.
\paragraph{Definition.}{\it A coproduct} on an $H$-algebra $A$
is a homomorphism of $H$-algebras
$\Delta: A\to A\wo A$.\\
\hbox to1em{\hfill}{\it A counit } on an $H$-algebra $A$ is
a homomorphism of $H$-algebras $\epe: A\to  D_T$.
\paragraph{}Finally, we can define the notions of an $H$-bialgebroid 
and an $H$-Hopf algebroid. 

\paragraph{Definition.}{\it An $H$-bialgebroid } is an
$H$-algebra $A$ equipped with a coassociative coproduct
$\Delta$ (i.e. such that $(\Delta\T \Id_A)\circ \Delta=
(\Id_A\T \Delta)\circ \Delta$, and a counit
$\epe$ such that $(\epe\T \Id_A)\circ \Delta=
(\Id_A\T \epe)\circ \Delta=\Id_A$.

\paragraph{}Let $A$ be an $H$-algebra. A linear map $S:A\to A$ is called 
{\it an antiautomorphism}
of $H$-algebras if it is an antiautomorphism
of algebras and $\mu_r\circ S=\mu_l,\, \mu_l\circ S=\mu_r$. From 
these conditions it follows that $S(A_{\al\beta})=A_{-\beta,-\alpha}$.

Let $A$ be an $H$-bialgebroid, and let $\Delta$, $\epe$
be the coproduct and counit of $A$. For $a\in A$, let
\begin{equation}\label{pres'}
\Delta(a)=\sum_i a^1_i\T a^2_i.
\end{equation}

\paragraph{Definition.} {\it{An antipode}} on the $H$-bialgebroid $A$
is an antiautomorphism of $H$-algebras $S:A\to A$ such that
for any $a\in A$ and any presentation \Ref{pres'} of $\Delta(a)$, one has
\begin{equation}
\sum_i a_i^1S(a_i^2)=\mu_l(\epe(a)1),\
\sum_i S(a_i^1)a_i^2=\mu_r(\epe(a)1),
\notag
\end{equation}
where $\epe(a)1\in M_{T}$ is the result of the application of the difference
operator $\epe(a)$ to the constant function $1$.

An $H$-bialgebroid with an antipode is called {\it an $H$-Hopf algebroid.}

\paragraph{Remarks.} 1. If $H=\C$ then the notions of $H$-algebra, $H$-bialgebroid, 
$H$-Hopf algebroid are the familiar notions of an algebra, bialgebra, and 
Hopf algebra. 

2. It is easy to see
that $D_T$ is an $H$-bialgebroid where $\Dl : D_T \to D_T\Tb D_T$
is the canonical isomorphism and $\epe = \Id$.
Furthermore, it is an $H$-Hopf algebroid with $S(D)=D^*$, where 
$D^*$ is the formal adjoint to the 
difference operator $D$ (i.e. $(f(\lambda)\mathcal T_\alpha)^*=
\mathcal T_\alpha^{-1}f(\lambda)$). This $H$-Hopf algebroid is an analog of 
the 1-dimensional Hopf algebra in the category of Hopf algebras. 

3. One can define the notions of an $H$-algebra, $H$-bialgebroid, 
$H$-Hopf algebroid if the group $T$ is not connected (for example, a finite 
group), in essentially the same way as above. More precisely, since 
in this case the algebra $M_T$ of meromorphic functions on $T$ 
is not a field but a direct sum of finitely many copies of a field,
one should introduce an additional axiom requiring that 
$A_{\alpha\beta}$ is a free module over $\mu_l(M_T)$ and $\mu_r(M_T)$.
Similarly, one can make all the above definitions in the case when $M_T$ 
is replaced with another algebra of functions on $T$ (rational functions, 
regular functions on some open set, etc.)

\paragraph{7.2. Dynamical representations of $H$-bialgebroids.}

One of the reasons $H$-bialgebroids are good analogs of bialgebras is 
that their representations, like representations of bialgebras, form 
a tensor category. However, these representations are not the usual 
representations but rather new objects which we call dynamical 
representations, and which we will now define. 

\paragraph{Definition.}{{\it A dynamical representation}} of an
$H$-algebra $A$ is a diagonalizable $H$-module
$W$ endowed with a homomorphism of $H$-algebras
$\pi_W:A\to D_{T, W}$, where $D_{T,W}$ is defined in Example 2. 

\paragraph{Definition.}{{\it A homomorphism}} of dynamical representations
$\phi: W_1\to W_2$ is an element of $\text{Hom}_{\C}(W_1,W_2\T M_{T})$
such that $\phi\circ \pi_{W_1}(x)=
\pi_{W_2}(x)\circ \phi$ for all $x\in A$.

\paragraph{Example.} If $A$ has a counit, then
$A$ has {\it the trivial representation}: $W=\C$, $\pi=\epe$.

\paragraph{}For diagonalizable $H$-modules $W,U$, 
let $f\in \text{Hom}(W,W\T M_{T})$ and \\
$g\in \text{Hom}(U, U\T
M_{T})$. Define $f\bo g \in \Hom (W\T U, W\T U \T M_{T})$ as
\begin{equation}\label{bar'}
f\bar{\T} g (\la)= f^{(1)}(\la-h^{(2)})(1\,\T g(\la))
\end{equation}
where $f^{(1)}(\la-h^{(2)})(1\,\T g(\la)) \, w\T u=
f(\la-\mu) w \,\T\,g(\la) u $ if $g(\la) u \in U[\mu]$.

\begin{lemma} [\cite{EV2}]  There is a natural embedding of
$H$-algebras\\
$\theta_{WU}:D_{T, W}\wo D_{T, U}\to D_{T, W\T U}$
(an isomorphism if $W,U$ are finite dimensional), given by the
formula $f\Ti_\beta\,\T g\,\Ti_\delta\to (f\bo g)\Ti_\delta$.
\end{lemma}

Now let us define the tensor product of dynamical representations
for $H$-bialgebroids.  
If $A$ is an $H$-bialgebroid, and 
$W$ and $U$ are two dynamical representations of $A$,
then we endow the $H$-module
$W\T U$ with the structure of a dynamical representation via
 $\pi_{W\T U}(x)=\theta_{WU}\circ (\pi_W\T \pi_U)\circ \Delta(x)$.
If $f:W_1\to W_2$ and $g:U_1\to U_2$
are homomorphisms of dynamical representations, then so is
$f\bo g :W_1\T U_1\to W_2\T U_2$.
Thus, dynamical representations of $A$ form a monoidal category
$\text{Rep}(A)$,
whose identity object is the trivial representation.

\paragraph{Remark.} 
If $A$ is an $H$-Hopf algebroid and $V$ is a dynamical representation,
then one can define the left and right dual dynamical representations
${}^*V$ and $V^*$. We will not discuss this notion here and refer 
the reader to \cite{EV2}. 

\paragraph{7.3. The $H$-bialgebroid associated to a function
$R:T\to \End(V\otimes~V)$.} Now, following \cite{EV2}, 
let us define an  $H$-bialgebroid $\bar A_R$ associated to a
meromorphic, zero weight function $R: T \to \End(V\T V)$,
where $V$ is a finite dimensional
diagonalizable
$H$-module (we assume that $R(\lambda)$ is nondegenerate for generic $\lambda$).
This is the dynamical analogue of the 
Faddeev-Reshetikhin-Takhtajan-Sklyanin construction of 
a bialgebra from an element $R\in \End(V\otimes~V)$, where 
$V$ is a vector space.

By definition, the algebra $\bar A_R$  is generated by two copies of
$M_{T}$ (embedded as subalgebras) and  matrix elements of the operator
$L\in \text{End}(V)\T \bar A_R$.
 We denote the elements of the first copy of
$M_{T}$ by $f(\la^1)$ and of the second copy by $f(\la^2)$,
where $f\in M_{T}$. We denote by $L_{\al\beta}$
the weight components of $L$ with respect to the natural
$T$-bigrading on $\End (V)$, so
that $L=(L_{\al\beta})$, where
$L_{\al\beta}\in \text{Hom}_\C(V[\beta],V[\al])\T \bar A_R$.

Introduce the moment maps for $\bar A_R$ by
$\mu_l(f)=f(\la^1)$, $\mu_r(f)=f(\la^2)$, and define the weight
decomposition by 
$$f(\la^1),f(\la^2)\in (\bar A_R)_{00},\qquad L_{\al\beta}\in
\text{Hom}_\C(V[\beta],V[\al])\T (\bar A_R)_{\al\beta}.$$

The defining relations for $\bar A_R$ are:
$$
f(\la^1)L_{\al\beta}=L_{\al\beta}f(\la^1+\al);\
f(\la^2)L_{\al\beta}=L_{\al\beta}f(\la^2+\beta);
[f(\la^1),g(\la^2)]=0;
$$
and the dynamical Yang-Baxter relation
\begin{align}\label{ybr}
R^{12}(\la^1)L^{13}L^{23}=:L^{23}L^{13}R^{12}(\la^2):.
\end{align}
Here the :: sign 
means that the matrix elements of $L$
should be put on the right of the matrix elements of $R$. Thus, if
$\{v_a\}$ is a homogeneous basis of $V$, and $L=\sum E_{ab}\T L_{ab}$,
$R(\la)(v_a\T v_b)=\sum R^{ab}_{cd}(\la)v_c\T v_d$, then \Ref{ybr} has the form
$$
\sum R^{xy}_{ac}(\la^1)L_{xb}L_{yd}=\sum R_{xy}^{bd}(\la^2)
L_{cy}L_{ax},
$$
where we sum over repeated indices.

Define the coproduct on $\bar A_R$, $\Delta: \bar A_R\to \bar A_R\Tb 
\bar A_R$, 
and the counit of $\bar A_R$ by
$$
\Delta(L)=L^{12}L^{13},\epe(L_{\alpha\beta})=\delta_{\alpha\beta}
Id_{V[\alpha]}\otimes \mathcal T_\alpha^{-1},
$$
where $\Id_{V[\al]}: V[\al]\to V[\al]$ is the identity operator.

\begin{prop}[\cite{EV2}]
$(\bar A_R,\Delta,\epe)$ is an $H$-bialgebroid.  
\end{prop}

\paragraph{Example.}Suppose that $R$ is the basic 
trigonometric solution of the QDYBE
(see Section 6). Then the defining relations for $\bar A_R$ look
like

\begin{align*}
f(\la^1)L_{bc}&=L_{bc}f(\la^1+\om_b),\\
f(\la^2)L_{bc}&=L_{bc}f(\la^2+\om_c), \\
L_{as}L_{at}&=\frac{\alpha_{st}(\la^2)}{1-\beta_{ts}(\la^2)}L_{at}L_{as},
s\ne t,\\
L_{bs}L_{as}&=\frac{\alpha_{ab}(\la_1)}{1-\beta_{ab}(\la^1)}L_{as}L_{bs},
a\ne b,\\
\al_{ab}(\la_1)L_{as}L_{bt}-\al_{st}(\la_2)L_{bt}L_{as}&=
(\beta_{ts}(\la_2)-\beta_{ab}(\la_1))L_{bs}L_{at}, a\ne b, s\ne t, 
\end{align*}
\noindent
where $\al_{ab}(\la)=\frac{q^{\la_a-\la_b+1}-1}{
q^{\la_a-\la_b}-1}$, $\beta_{ab}(\la)=\frac{q-1}{q^{\la_a-\la_b}-1}$. 

\paragraph{}We note that we don't need any special properties of $R$
(like the dynamical Yang-Baxter equation or Hecke condition) to define 
the $H$-bialgebroid $\bar A_R$. However, if we take a "randomly chosen" 
function $R$, the $H$-bialgebroid $\bar A_R$ will most likely have rather bad 
properties; i.e. it will be rather small and will not have interesting 
dynamical representations. The simplest way to ensure the existence 
of at least one interesting dynamical representation is to require that 
$R$ satisfies the QDYBE. This
is so because  of the following proposition. 
\paragraph{}If $(W, \pi_W)$ is a dynamical representation of an $H$-algebra $A$,
we denote  $\pi_W^0:A\to \text{Hom}(W,W\T M_T)$
the map defined by $\pi_W^0(x)w=
\pi_W(x)w$, $w\in W$
 (the difference operator $\pi_W(x)$ restricted to the constant functions).
It is clear that $\pi_W$ is completely determined by $\pi_W^0$.
\begin{prop} If $R$ satisfies the QDYBE then $\bar
A_R$ has a dynamical representation realized in the space $V$, 
with
$\pi_V^0(\lambda)=R(\lambda)$. 
\end{prop}

This representation is called {\it the vector representation}. 

However, even if $R$ satisfies the QDYBE,  the
$H$-bialgebroid $\bar A_R$ may not be completely
satisfactory.  In particular, one may ask the following question:
does $\bar A_R$ define  a "good quantum matrix algebra"? More
precisely, does the  Hilbert series of
$\bar A_R$ equal to that of the function ring on the matrix
algebra, i.e.
$(1-t)^{-dim(V)^2}$? In general, the answer is no, even if the
quantum dynamical Yang-Baxter  equation is satisfied. 

In fact, here is the place where the Hecke condition 
comes handy. Namely, we have the following proposition, which 
is a generalization of a well known proposition in the theory 
of quantum groups (due to Faddeev, Reshetikhin, Takhtajan). 

\begin{prop}[\cite{EV2}] Suppose that $R$ satisfies the QDYBE
and the Hecke condition with $q$ not equal to a  nontrivial root
of $1$. Then  the space $\bar A_R^m$ of polynomials of degree
$m$ in generators 
$L_{\alpha\beta}$ in $\bar A_R$ is a free $M_T\otimes M_T$-module, and 
the ranks of these modules are given by
$$
\sum_{m\ge 0} rk(\bar A_R^m)t^m=(1-t)^{-n^2}.
$$
\end{prop}
\paragraph{Remark.}The proof of this proposition, like the proof of its nondynamical analog, 
 is based on the fact that under the assumptions 
of the Proposition, $\bar A_R^m$ is a representation of the Hecke algebra.
This justifies the name "Hecke condition". 

\paragraph{7.4. The $H$-Hopf algebroid $A_R$.}Suppose now that $R$ satisfies the
QDYBE and the Hecke condition where $q$ is not a nontrivial root of
$1$. In this case, it turns out that, analogously to  the nondynamical
case, a suitable localization
$A_R$ of $\bar A_R$  is actually a Hopf algebroid. Namely,
define $A_R$ by adjoining  to $\bar A_R$ a new element
$L^{-1}$, with the relation
$LL^{-1}=L^{-1}L=1$. It is easy to see that the structure 
of an $H$-bialgebroid on $\bar A_R$ naturally extends to 
$A_R$, and it can be shown that $A_R$ admits a unique antipode $S$ such that 
$S(L)=L^{-1}$. This antipode equips $A_R$ with a structure of an 
$H$-Hopf algebroid. This $H$-Hopf algebroid is a quantization 
of the group $GL_n$ in the same sense as the $H$-bialgebroid 
$\bar A_R$ is a quantization of the matrix algebra $Mat_n$.  

\paragraph{7.5. Quasiclassical limit.} In conclusion of the section, we would like to explain 
why the $H$-Hopf algebroid $A_R$ considered here 
(for the basic rational or trigonometric solution $R$ of the 
QDYBE) should be  
regarded as a quantization of the Poisson groupoid $X_r$ corresponding 
to the basic rational, respectively trigonometric, solution $r$ 
of the CDYBE (for the 
definition of $X_r$, see Section 5). 

To see this, consider the $H$-Hopf algebroid $A_R$ 
with $T=\h^*$ for some finite dimensional abelian Lie algebra $\h$, and 
$M_T$ replaced with the ring of regular functions on some open subset 
$U$ of $\h^*$. Introduce a formal parameter $\gamma$ (like in Section 3), 
and make a change of variable $\la\to \la/\gamma$ in the defining 
relations for $A_R$. It is easy to see that 
the resulting algebra $A_R^\gamma$ 
over $\C[[\gamma]]$ is a deformation of a commutative algebra. 
The above result about the Hilbert series implies that 
this deformation is flat, so the quotient algebra $A_R^0:=A_R^\gamma/(\gamma)$ 
obtains a Poisson structure. Let $X$ be the spectrum of $A_R^0$;
it is an algebraic Poisson manifold. It is not difficult to show that the 
the moment maps, coproduct, 
counit, and antipode of $A_R$ define maps $s,t,m,E,i$ (see Section 5) 
for $X$, which equips $X$ with the structure of a Poisson groupoid with base 
$U$. 
Moreover, it is easy to check that the Poisson groupoid 
$X$ is naturally isomorphic 
to $X_r$.   

\section{The 
universal fusion matrix and the Arnaudon-Buffenoir-Ragoucy-Roche equation}
\paragraph{8.1. The ABRR equation.} In \cite{ABRR}, Arnaudon, Buffenoir, Ragoucy and Roche give a general method for constructing the universal fusion matrix $J(\lambda)$, which lives in some completion of $U_q(\g)^{\otimes 2}$, 
i.e the unique element satisfying $J_{VW}(\lambda)=J(\lambda)_{|V\otimes W}$ for all $V,W$. A similar approach is suggested in \cite{JKOS}, based on 
the method of \cite{Fr}\\
\hbox to1em{\hfill}Let $U'(\b_{\pm})$ be the kernel of the projection
$U(\b_\pm) \to U({\h})$. We use the same notations with
the index $q$ for the quantum analogs of these objects. We set
$\theta(\lambda)=\lambda +
\rho-\frac{1}{2}\sum_{i}x_i^2 \in U\h$ where as usual
$\rho=\frac{1}{2}\sum_{\alpha \in \Delta^+}h_\alpha$ and $(x_i)$ is
an orthonormal basis of $\h$. Set
$\mathcal{R}_0=\mathcal{R}q^{-\sum x_i \otimes x_i}$. It is known
that $\Ro \in 1 + U'_q(\b_+) \otimes U'_q(\b_-)$.
\begin{theo}[\cite{ABRR}] The universal fusion matrix $J(\lambda)$
of $U_q(\g)$ is the unique solution of the form $1+ U'_q(\b_-) \otimes
U'_q(\b_+)$ of the equation
\begin{equation}\label{E:abrrq}
J(\lambda)(1 \otimes q^{2\theta(\lambda)})=\Ro^{21}(1 \otimes
q^{2\theta(\lambda)})J(\lambda).
\end{equation}
 The universal fusion matrix $J(\lambda)$ of $U(\g)$ is the unique solution of the form $1+
U'(\b_-) \otimes U'(\b_+)$ of the equation
\begin{equation}\label{E:abrr}
[J(\lambda),1 \otimes \theta(\lambda)]=(\sum_{\alpha \in \Delta^+} e_{-\alpha}\otimes e_{\alpha}) J(\lambda)
\end{equation}
\end{theo}
\noindent

We will call these equations the
\textit{ABRR equations} for $U_q(\g)$ and for $\g$, respectively.

\textit{Proof.} Let us first show the statement about uniqueness.
Let $T(\lambda)\in 1 + U'_q(\b_-) \otimes U'_q(\b_+)$ be any
solution of (\ref{E:abrrq}). Then
\begin{align*}
(\Ro^{21})^{-1}T(\lambda)= \mathrm{Ad}\;(1 \otimes q^{2\theta(\lambda)})T(\lambda)\\
\Leftrightarrow ((\Ro^{21})^{-1}-1)T(\lambda)&= ( \mathrm{Ad}\;(1 \otimes q^{2\theta(\lambda)})-1)T(\lambda)\\
\Leftrightarrow T(\lambda)=1 + ( \mathrm{Ad}\;(1 &\otimes  q^{2\theta(\lambda)})-1)^{-1}((\Ro^{21})^{-1}-1)T(\lambda)
\end{align*}
Now notice that $((\Ro^{21})^{-1}-1) \in U'_q(\b_-) \otimes
U'_q(\b_+)$. This implies that $T(\lambda)$ can be recusively
constructed as follows. Set $T_0(\lambda)=1$ and put
$$T_{n+1}(\lambda)=1 + ( \mathrm{Ad}\;(1 \otimes
q^{2\theta(\lambda)})-1)^{-1}((\Ro^{21})^{-1}-1)T_n(\lambda).$$
Then $\mathrm{lim}_{n \to \infty} T_n(\lambda)=T(\lambda)$ (the
limit is in the sense of stabilization). In particular there exists a
unique solution to (\ref{E:abrrq}) of the given form.\\
\hbox to1em{\hfill}The proof in the rational case (i.e in the case of
a simple Lie algebra $\g$) is similar. In that case, the recursive
construction is given by $T_0(\lambda)=1$ and
$$T_{n+1}(\lambda)=1-\mathrm{ad}(1 \otimes
\theta(\lambda)^{-1})(\sum_{\alpha \in \Delta^+} e_{-\alpha}
\otimes e_{\alpha})T_{n}(\lambda).$$ We now give a proof that the
fusion matrix $J(\lambda)$ actually satisfies the ABRR relation in
the case of simple Lie algebras. The proof in the case of quantum
groups is analogous but technically more challenging, and is given  in
Appendix B. 

Let $C$ be the quadratic Casimir operator in the center of the universal enveloping algebra $U\g$:
$$C=\sum_i x_i^2 + 2\rho + 2 \sum_{\alpha \in \Delta^+} e_{-\alpha}e_\alpha.$$
Then $C$ acts on any highest weight representation of $\g$ of highest weight $\lambda$ by the scalar $(\lambda,\lambda+2\rho)$. Now let $V,W$ be two finite-dimensional $\g$-modules and let $v \in V,\;w \in W$ be two homogeneous elements of weight $\mathrm{wt}(v)$ and $\mathrm{wt}(w)$. We compute the quantity
$$F(\lambda)=\langle x_{\lambda-\mathrm{wt}(v)-\mathrm{wt}(w)}^*, \Phi^w_{\lambda-\mathrm{wt}(v)} (C \otimes 1) \Phi^v_\lambda x_\lambda \rangle$$
in two different ways. On one hand we have
\begin{equation}\label{E:61}
F=(\lambda-\mathrm{wt}(v),\lambda-\mathrm{wt}(v)+2\rho) J(\lambda)(w \otimes v).
\end{equation}
On the other hand,  
\begin{equation*}
\begin{split}
F=\langle x_{\lambda-\mathrm{wt}(v)-\mathrm{wt}(w)}^*,& \big\{2\big(\sum_\alpha((e_{-\alpha} e_\alpha)_{1} +  (e_{-\alpha} e_\alpha)_{2} + (e_\alpha \otimes e_{-\alpha})_{12} + (e_{-\alpha} \otimes e_\alpha)_{12}) + \rho_1
+\rho_2)\big) \\ 
&+\sum_i (x_i^2)_1+(x_i^2)_2+2(x_i \otimes x_i)_{12} \big\} \Phi^w_{\lambda-\mathrm{wt}(v)}
\Phi^v_\lambda x_\lambda \rangle.
\end{split}
\end{equation*}
Since $ x_{\lambda-\mathrm{wt}(v)-\mathrm{wt}(w)}^*$ is a lowest weight vector, it is clear that
$$\langle x_{\lambda-\mathrm{wt}(v)-\mathrm{wt}(w)}^*, (e_{-\alpha} e_\alpha)_1 \Phi^w_{\lambda-\mathrm{wt}(v)} \Phi^v_\lambda x_\lambda \rangle=\langle x_{\lambda-\mathrm{wt}(v)-\mathrm{wt}(w)}^*, (e_{-\alpha} \otimes e_\alpha)_{12} \Phi^w_{\lambda-\mathrm{wt}(v)} \Phi^v_\lambda x_\lambda \rangle=0.$$
Moreover, by the intertwining property again, we have
\begin{align*}
(e_\alpha \otimes e_{-\alpha})_{12} \Phi^w_{\lambda-\mathrm{wt}(v)} \Phi^v_\lambda x_\lambda &= (-(e_{-\alpha} e_\alpha)_2 - ( e_{-\alpha} \otimes e_\alpha)_{23}) \Phi^w_{\lambda-\mathrm{wt}(v)} \Phi^v_\lambda
x_\lambda,\\
 (\rho_1 + \rho_2) \Phi^w_{\lambda-\mathrm{wt}(v)} \Phi^v_\lambda x_\lambda &= -\rho_3 
\Phi^w_{\lambda-\mathrm{wt}(v)} \Phi^v_\lambda v_\lambda + \Phi^w_{\lambda-\mathrm{wt}(v)}
\Phi^v_\lambda \rho x_\lambda, 
\end{align*}
and
\begin{equation*}
\begin{split} 
\sum_i((x_i^2)_1 +&( x_i^2)_2 +2( x_i \otimes x_i)_{12}) 
\Phi^w_{\lambda-\mathrm{wt}(v)} \Phi^v_\lambda x_\lambda \\
&=-\sum_i\big(2(x_i \otimes x_i)_{13} +2(x_i \otimes
x_i)_{23} + (x_i^2)_3) \Phi^w_{\lambda-\mathrm{wt}(v)} \Phi^v_\lambda x_\lambda+\Phi^w_{\lambda-\mathrm{wt}(v)} \Phi^v_\lambda x_i^2 x_\lambda\big) \\
&=\sum_i\big(\Phi^w_{\lambda-\mathrm{wt}(v)}
\Phi^v_\lambda x_i^2 x_\lambda + ((x_i^2)_3 - 2\lambda_3)\Phi^w_{\lambda-\mathrm{wt}(v)}
\Phi^v_\lambda x_\lambda\big).
\end{split}
\end{equation*}
Summing up these equations, we finally obtain
\begin{equation}\label{E:62}
F=(-2\rho_{|V} + 2(\rho,\lambda) + (\sum_i x_i^2)_{|V} -2\lambda_{|V} + \lambda^2-2\sum_{\alpha \in \Delta^+}e_{-\alpha} \otimes e_\alpha)J(\lambda)(w \otimes v)
\end{equation}
Combining (\ref{E:61}) and (\ref{E:62}) yields
$$(\sum_i x_i^2-2(\lambda+\rho))_{|V}J(\lambda) - (\mathrm{wt}(v)^2-2(\lambda+\rho,\mathrm{wt}(v)))J(\lambda)=-2(\sum_{\alpha \in \Delta^+}e_{-\alpha} \otimes e_\alpha)J(\lambda)$$
which is equivalent to (\ref{E:abrr}).\qed
\paragraph{Example.} Let us use the recursive procedure in the proof above 
to compute $J(\lambda)$ for $U\mathfrak({sl}_2)$. Setting $J(\lambda)=1 + \sum_{n \geq 1} J^{(n)}(\lambda)$ where $J^{(n)}(\lambda)\in
U'(\b_-)[-2n] \otimes U'(\b_+)[2n]$, the ABRR equation reads
$$\frac{1}{2}[1 \otimes ((\lambda +1)h-\frac{h^2}{2}),J-1]=-(f \otimes e)J,$$
which gives the recurrence relation
$$1 \otimes ((\lambda+1)n-nh + n^2)J^{(n)}=(-f \otimes e)J^{(n-1)}.$$
Hence
$$J^{(n)}(\lambda)=\frac{(-1)^n}{n!}f^n \otimes (\lambda-h+n+1)^{-1} \ldots (\lambda-h +2n)^{-1}e^n.$$
This formula and its quantum analogue 
were obtained in the pioneering paper \cite{BBB}, which was a motivation 
to the authors of \cite{ABRR}.

\paragraph{8.2. Classical limits of fusion and exchange matrices.} The ABRR relations can also be used to derive the classical limits of the fusion (and thus of the exchange) matrices. Setting $q=e^{-\gamma/2}$, rescaling $\lambda \mapsto \frac{\lambda}{\gamma}$ and considering the limit $\gamma \to 0$ yields the following classical version of the ABRR equation:
$$\mathrm{Ad}( 1 \otimes e^{-\lambda})j(\lambda)-j(\lambda)=\sum_{\alpha \in \Delta^+} e_{-\alpha} \otimes e_{\alpha},$$
which admits the unique lower triangular solution
$$j(\lambda)=-\sum_{\alpha \in \Delta^+}\frac{e_{-\alpha} \otimes e_\alpha}{1-e^{-(\alpha,\lambda)}}.$$
 From this we deduce
$$r(\lambda)=r^{21}+j(\lambda)-j^{21}(\lambda)=\frac{1}{2}\Omega + \frac{1}{2}\sum_{\alpha >0}\mathrm{cotanh} \left( \frac{1}{2} \langle \alpha,\lambda\rangle \right) (e_\alpha \otimes e_{-_\alpha} - e_{-\alpha} \otimes e_\alpha)$$
which is consistent with Proposition 3.3.
\paragraph{}The case of a simple Lie algebra $\g$ is completely analogous; the classical version of the ABRR equation is 
$$[j(\lambda),\lambda \otimes 1]=\sum_{\alpha \in \Delta^+} e_{-\alpha} \otimes e_{\alpha},$$
which admits the unique solution
$$j(\lambda)=-\sum_{\alpha \in \Delta^+}\frac{e_{-\alpha} \otimes e_\alpha}{(\lambda,\alpha)}.$$
This yields Proposition 3.2. 

\section{Transfer matrices and Generalized Macdonald-Ruijsenaars equations}
\paragraph{9.1. Transfer matrices.} We first recall the well-known transfer matrix construction. 

Let $A$ be a Hopf algebra with a commutative Grothendieck ring,
and let ${\mathcal R} \in A \otimes A$ be 
an element such that $(\Delta\otimes 1)({\mathcal R})={\mathcal R}^{13}{\mathcal R}^{23}$. A basic example is: 
$A$ is quasitriangular, ${\mathcal R}$ is its universal R-matrix. 

For any finite-dimensional representation $\pi_V: A \to \mathrm{End}\;V$ of $A$, set
$$T_V=\mathrm{Tr}_{|V} (\pi_V \otimes 1)({\mathcal R}) \in A.$$
These elements are called transfer matrices. 

\begin{lem}\label{L:tran} For any finite-dimensional $A$-modules $V,W$ we have $T_VT_W=T_{V
\otimes W}=T_WT_V$. 
 \end{lem}
\noindent
\textit{Proof.} By definition we have 
$$(\pi_{V \otimes W} \otimes 1){\mathcal R}=(\pi_V \otimes \pi_W \otimes 1)(\Delta \otimes 1){\mathcal R}=(\pi_V \otimes 1){\mathcal R}_{13}(\pi_W \otimes 1){\mathcal R}_{23},$$
which implies the first equality. The second equality follows 
from the commutativity of the Grothendieck ring.  
\qed
\paragraph{}The transfer matrix construction gives rise to interesting examples of 
quantum integrable systems which arise in statistical mechanics. 
For example, if $A$ is the quantum affine algebra or the elliptic algebra, 
one gets transfer matrices of the 6-vertex and 8-vertex models, respectively.

\paragraph{}We adapt the notion of transfer matrices in our dynamical setting in the following way. Let $\g$ be a simple Lie algebra and let $U_q(\g)$ be the associated quantum group. For any two finite-dimensional $U_q(\g)$-modules $V$ and $W$ let $R_{VW}(\lambda)$ be the exchange matrix. It is more convenient to work with the shifted exchange matrix $\R(\lambda)=R(-\lambda-\rho)$.\\
\hbox to1em{\hfill}Let $\mathcal{F}_V$ be the space of
$V[0]$-valued meromorphic functions on $\h^*$. For $\nu
\in \h^*$ let $T_\nu \in \mathrm{End}(\mathcal{F}_V)$ be the shift operator $(T_\nu
f)(\lambda)=f(\lambda+\nu)$.  As pointed out in \cite{FV3},  the role of the transfer matrix is played by the
following difference operator
$$\mathcal D^V_W=\sum_\nu \mathrm{Tr}_{|W[\nu]}(\R_{WV}(\lambda))T_\nu.$$
It follows from the dynamical 2-cocycle condition 
for fusion matrices (see Proposition 2.3) that for any
$U_q(\g)$-modules $U,V,W$ we have
$$\mathcal D^U_{V \otimes W}=\mathcal D^U_V \mathcal D^U_W=\mathcal D^U_W \mathcal D^U_V.$$
Hence $\{\mathcal D^U_W\}$ span a commuting family of difference operators acting on $\mathcal{F}_U$.
\paragraph{9.2. Weighted trace functions.} Let $V$ be a finite-dimensional $U_q(\g)$-module. Recall that, for any homogeneous vector $v \in V[\nu]$ and for generic $\mu \in \h^*$ there exists a unique intertwiner $\Phi^v_\mu: M_\mu \to M_{\mu-\nu} \otimes V$ such that $\langle v_{\mu-\nu}^*, \Phi^v_\mu v_\mu\rangle=v$. Set
$$\Phi^V_\mu=\sum_{v \in \mathcal{B}} \Phi^v_\mu \otimes v^* \in \mathrm{Hom}_\C\big( M_\mu, \bigoplus_\nu M_{\mu-\nu} \otimes V[\nu] \otimes V^*[-\nu]\big),$$
where $\mathcal{B}$ is any homogeneous basis of $V$. Consider the weighted trace function
$$\Psi_V(\lambda,\mu)=\mathrm{Tr}\;(\Phi^V_\mu q^{2\lambda}) \in V[0] \otimes V^*[0]$$
where $q^{2\lambda}$ acts on any $\h$-semisimple $U_q(\g)$-module $U$ by
$q^{2\lambda}_{|U[\nu]}=q^{2(\lambda,\nu)}Id$. It can be shown that $\Psi_V \in V[0] \otimes V^*[0] \otimes
q^{2(\lambda,\mu)} \C(q^\lambda) \otimes \C(q^\mu)$. Let 
$$\delta_q(\lambda)=\big(\mathrm{Tr}_{|M_{-\rho}}(q^{2\lambda})\big)^{-1}=q^{-2(\lambda,\rho)} \prod_{\alpha >0} (1-q^{-2(\lambda,\alpha)})$$
be the Weyl denominator, and set
$$Q(\lambda)=m^{op}(1 \otimes S^{-1})(J(-\lambda-\rho)),$$
where $m^{op}: U_q(\g) \otimes U_q(\g) \to U_q(\g), a \otimes b \mapsto ba$ and where $J(\lambda)$ is the universal fusion matrix. It can be shown that $Q(\lambda)$ is invertible. Finally, set
$$F_V(\lambda,\mu)=Q^{-1}(\mu)_{|V^*} \Psi_V(\lambda,-\mu-\rho)
\delta_q(\lambda).$$
\begin{theo} [\cite{EV4}, 
The Macdonald-Ruijsenaars equations]\label{mr}
For any two finite-dimensional $U_q(\g)$-modules, we have
$$
\mathcal D^{\lambda,V}_WF_V(\lambda,\mu)= \chi_W(q^{-2\mu})F_V(\lambda,\mu)
$$
where $\chi_W(q^x)=\sum\;\mathrm{dim}\;W[\nu] q^{(\nu,x)}$ is the character of $W$.
\end{theo}
\begin{theo}[\cite{EV4}, The dual Macdonald-Ruijsenaars equations]\label{dmr}
For any two finite-dimensional $U_q(\g)$-modules, we have
$$
\mathcal D^{\mu,V^*}_WF_{V}(\lambda,\mu)= \chi_W(q^{-2\lambda})F_V(\lambda,\mu) 
$$
\end{theo}
\noindent
In the above, we add a superscript to $\mathcal D$ to specify on which variable the difference operators act. Thus, in Theorem \ref{mr}, $\mathcal D^V_W$ acts on functions of the variable $\lambda$ in the component $V[0]$, and in Theorem 
\ref{dmr}, $\mathcal D^{V^*}_W$ acts on functions in the variable $\mu$ in the component $V^*[0]$.

 From Theorems 9.1 and 9.2 it is not difficult to deduce the following result:

\begin{theo}[\cite{EV4}, The symmetry identity]\label{sym}
For any finite-dimensional $U_q(\g)$-module we have
$$F_V(\lambda,\mu)=F^*_{V^*}(\mu,\lambda),$$
where $*: V[0] \otimes V^*[0] \to V^*[0] \otimes V[0]$ is the permutation.
\end{theo}

\paragraph{9.3. Relation to Macdonald theory.}Let us now restrict ourselves to the case of $\g=sl_n$,
 and let $V$ be the q-analogue of the representation $S^{mn}\C^n$.
The zero-weight subspace of this representation is 1-dimensional, 
so the function $\Psi_V$ can be regarded as a scalar function. 
We will denote this scalar function by $\Psi_m(q,\la,\mu)$.

Recall the definition of Macdonald operators \cite{Ma,EK1}. 
They are operators on the space of functions 
$f(\la_1,...,\la_n)$ which are invariant under simultaneous shifting 
of the variables, $\la_i\to \la_i+c$, and have the form
$$
M_r=\sum_{I\subset \{1,...,n\}: |I|=r}
\left(\prod_{i\in I,j\notin I}
\frac{tq^{2\la_i}-t^{-1}q^{2\la_j}}{q^{2\la_i}-q^{2\la_j}}\right)T_I,
$$  
where $T_I\la_j=\la_j$ if $j\notin I$ and $T_I\la_j=\la_j+1$ if 
$j\in I$. Here $q,t$ are parameters. We will assume that $t=q^{m+1}$, 
where $m$ is a nonnegative integer. 

It is known \cite{Ma} that the operators $M_r$ commute. 
 From this it can be deduced that 
for a generic $\mu=(\mu_1,...,\mu_n)$, $\sum \mu_i=0$, there exists 
a unique power series $f_{m0}(q,\la,\mu)\in \C[[q^{\la_2-\la_1},...,
q^{\la_{n}-\la_{n-1}}]]$ such that the series
$f_m(q,\la,\mu):=q^{2(\la,\mu-m\rho)}f_{m0}(q,\la,\mu)$
satisfies difference equations 
$$
M_rf_m(q,\la,\mu)=
(\sum_{I\subset \{1,...,n\}: |I|=r} q^{2\sum_{i\in I}(\mu+\rho)_i})
f_m(q,\la,\mu).
$$

\paragraph{Remark.}The series $f_{m0}$ is convergent to an analytic (in fact, a trigonometric)
function.

The following theorem is contained in \cite{EK1}.

\begin{theo}[\cite{EK1}, Theorem 5]\label{trace} One has 
$$
f_m(q,\la,\mu)=\gamma_m(q,\la)^{-1}\Psi_m(q^{-1},-\la,\mu),
$$
where 
$$
\gamma_m(q,\la):=
\prod_{i=1}^m\prod_{l<j}(q^{\la_l-\la_j}-q^{2i}q^{\la_j-\la_l}).
$$
\end{theo}

Let $\mathcal D_W(q^{-1},-\la)$ denote the difference operator, obtained 
from the operator $\mathcal D_W$ defined in Section 1 by the transformation 
$q\to q^{-1}$ and the change of coordinates $\la\to -\la$.
Let $\Lambda^r\C^n$ denote the q-analog of the r-th fundamental representation 
of $sl_n$. 
 
\begin{cor}
$$
\mathcal D_{\Lambda^r\C^n}(q^{-1},-\la)=\delta_q(\la)\gamma_m(q,\la)\circ M_r
\circ \gamma_m(q,\la)^{-1}\delta_q(\la)^{-1}.
$$
\end{cor}
\noindent
\textit{Proof.} This follows from Theorem \ref{trace} and Theorem \ref{mr}. 

\paragraph{Remark.} In the theorems of this section, Verma modules 
$M_\mu$ can be replaced with 
finite dimensional irreducible modules $L_\mu$ with sufficiently 
large highest weight, and one can prove analogs of these theorems
in this situation (in the same way as for Verma modules). 
In particular, one may
set $\hat\Psi_m(q,\la,\mu)
=\text{Tr}(\hat\Phi_\mu^V q^{2\la})$, where 
$\hat \Phi_\mu^V:L_\mu\to L_\mu\otimes V\otimes V^*[0]$ is the intertwiner with 
highest coefficient 1 (Such an operator exists iff $\mu-m\rho\ge 0$, 
see \cite{EK1}). Then one can show analogously to Theorem 9.1 
(see \cite{EK1}) that 
the function 
$\hat f_m(q,\la,\mu):=\gamma_m(q,\la)^{-1}\hat\Psi_m(q^{-1},-\la,\mu+m\rho)$ 
is the Macdonald polynomial $P_\mu(q,t,q^{2\la})$ with highest weight $\mu$
($\mu$ is a dominant integral weight). 
In this case, Theorem \ref{mr} says that Macdonald's polynomials 
are eigenfunctions of Macdonald's operators, Theorem \ref{dmr} gives recursive 
relations for Macdonald's polynomials with respect to the weight 
(for $sl(2)$ -- the usual 3-term relation for orthogonal polynomials), 
and Theorem \ref{sym} is the Macdonald symmetry identity (see \cite{Ma}).

\section{Appendix A: Classical dynamical r-matrices on a simple Lie
algebra $\g$ with respect to\\ $\l \subset \h$.} 
\paragraph{A.1.} Let $\g$ be a simple
complex Lie algebra and $\l$ a commutative subalgebra of $\g$
consisting of semisimple elements. Then $\l \subset \h$ for some
Cartan subalgebra $\h$. We keep the notations of Section 1. In this
appendix, we give a classification of all classical dynamical
r-matrices $\l^* \to (\g \otimes \g)^\l$ with coupling constant
$1$. Note that we can suppose without loss of generality that the
restriction of $\langle\,,\,\rangle$ to $\l$ is nondegenerate. Indeed,
given a dynamical r-matrix $r: \l^* \to (\g \otimes \g)^\l$, we can
always replace $\l$ by the largest subalgebra of $\h$ under which
$r$ is invariant, and this subalgebra is real. 

\paragraph{A.2. Gauge transformations.} Let $\Omega' \in \l^\perp
\otimes \l^\perp$ (where the orthogonal complement is in $\g$) be
the inverse (Casimir) element to the form $\langle\;,\;\rangle$. If
$r(\lambda)=\frac{1}{2}\Omega + (\varphi(\lambda) \otimes
1)\Omega'$, $\varphi: \l^* \to \mathrm{End}_\l(\l^\perp)$ is a
meromorphic function with values in $(\g \otimes \g)^\l$, and if $f:
\l^* \to \h$ is any meromorphic function, set
$$r^f(\lambda)=\frac{1}{2}\Omega + (e^{-ad \;f(\lambda)}\varphi(\lambda)e^{ad\;f(\lambda)} \otimes 1)\Omega'.$$
\begin{lem} The transformations $r(\lambda) \mapsto r^f(\lambda)$ preserve the set of classical dynamical r-matrices with coupling constant $1$.
\end{lem}
Two r-matrices which can be obtained one from the other by such a transformation are called
\textit{gauge-equivalent}.
\paragraph{A.3. Classification of dynamical r-matrices.} Let $\h$ be a Cartan subalgebra of $\g$, and let $\Pi
\subset \h^*$ be a system of simple roots in $\Delta$. Let
$\h_0\subset \h$ be the orthogonal complement of $\l$ in $\h$.
\paragraph{Definition.} A \textit{generalized Belavin-Drinfeld triple} is a triple $(\Gamma_1,\Gamma_2,\tau)$ where $\Gamma_1, \Gamma_2 \subset \Pi$, and where $\tau: \Gamma_1 \stackrel{\sim}{\to} \Gamma_2$ is a norm-preserving isomorphism.\\
\paragraph{}Given a generalized Belavin-Drinfeld triple
$(\Gamma_1,\Gamma_2,\tau)$, we extend linearly the map $\tau$ to a
norm-preserving bijection $\langle \Gamma_1 \rangle \to \langle
\Gamma_2\rangle$, where $\langle \Gamma_1 \rangle$ (resp. $\langle
\Gamma_2\rangle$) is the set of roots $\alpha \in \Delta$ which are
linear combinations of simple roots from $\Gamma_1$ (resp. from $\Gamma_2$).\\
\hbox to1em{\hfill}We will say that a generalized Belavin-Drinfeld
triple $(\Gamma_1,\Gamma_2,\tau)$ is $\l$-admissible if
$\tau(\alpha)-\alpha\in \l^\perp$ for all $\alpha \in \Gamma_1$, if
$\tau$ satisfies
the following condition:
for every cycle $\alpha \mapsto \tau(\alpha) \mapsto \ldots
\mapsto \tau^r(\alpha)=\alpha$ we have
$\alpha+\tau(\alpha)+\ldots +\tau^{r-1}\alpha \in \l$.
\paragraph{}If $(\Gamma_1,\Gamma_2,\tau)$ is a generalized Belavin-Drinfeld triple, let $\g_{\Gamma_1},\;\g_{\Gamma_2}$ be the subalgebras generated by
$e_\alpha,f_\alpha\;\alpha \in\Gamma_1$ (resp. generated by $e_\alpha,f_\alpha\;\alpha \in\Gamma_2$).
The map
$\tau: e_\alpha \mapsto e_{\tau(\alpha)},\;\alpha \in \Gamma_1$ extends to an isomorphism
$\tau:\;\g_{\Gamma_1} \stackrel{\sim}{\to} \g_{\Gamma_2}$. Finally, define an operator $K: \l^* \to
\mathrm{Hom}\,(\g_{\Gamma_1},\g)$ by
$$K(\lambda)e_{\alpha}=\sum_{n>0}e^{-n(\alpha,\lambda)}\tau^n(e_\alpha).$$
Notice that this sum is finite if $\tau$ acts nilpotently on $\alpha$.
\begin{theo}[\cite{S}]\label{T:app} Let $(\Gamma_1,\Gamma_2,\tau)$ be
an $\l$-admissible generalized Belavin-Drinfeld triple. 

(i) The equation 
$$((\alpha-\tau(\alpha)) \otimes 1 )r_0=\frac{1}{2}((\tau(\alpha)+\alpha)
\otimes 1)
\Omega_{\h_0},$$  where $\Omega_{\h_0}\subset \h_0 \otimes
\h_0$ is the inverse element to the form $\langle\,,\,\rangle$, has solutions
$r_0\in \Lambda^2\h_0$.  

(ii) Let $r_0 \in
\Lambda^2\h_0$ satisfy the equation from (i). 
Then
\begin{equation}\label{E:app}
r(\lambda)=\frac{1}{2}\Omega + r_0 + \sum_{\underset{e_\alpha \in \g_{\Gamma_1}}{\alpha\in \Delta^+}} K(\lambda)e_\alpha \wedge
f_{\alpha} + \sum_{\alpha\in \Delta^+}
\frac{1}{2}e_\alpha
\wedge f_{\alpha}
\end{equation}
is a classical dynamical r-matrix. Conversely, any classical dynamical
r-matrix $r:\l^* \to (\g \otimes \g)^\l$ with coupling constant $1$ is
gauge equivalent to one of the above form, for suitable choices of
Cartan subalgebra $\h$ containing $\l$, polarization of $\g$ and
$\l$-admissible generalized Belavin-Drinfeld triple.
\end{theo}

{\textit{Proof.}} Let us prove statement (i); statement (ii) is proved in 
\cite{S}. 
Let $\l_{max}$ be the Lie algebra of all 
$x\in \h$ such that $(\al-\tau(\al),x)=0$, 
and let ${\frak p}$ be the orthogonal complement of $\l$ in $\l_{max}$. 
Then we have an orthogonal direct sum decomposition
$\h_0={\frak p}\oplus\l_{max}^\perp$. Let us regard $r_0$ as a bilinear 
form on $\h$ (via the standard inner product). The equation from (i) determines 
$r_0(x,y)$ where $x\in \l_{max}^{\perp}$, and $y$ is arbitrary.
To check that $r_0$ can be extended to a skew symmetric form, it suffices 
to check that it is skew-symmetric on $\l_{max}^\perp$. 
But using the equation from (i) we find that 
$$
r_0(\al-\tau(\al),\beta-\tau(\beta))=1/2(\al+\tau(\al),\beta-\tau(\beta))=
\frac{1}{2}((\beta,\tau(\al))-(\al,\tau(\beta)))
$$
(we use that $\tau$ preseves $(,)$), which is obviously skew symmetric. 
\qed

\paragraph{Remarks.} 1. This classification is very similar in spirit to the classification of classical r-matrices $r \in \g \otimes \g$ satisfying $r+r^{21}=\Omega$ (quasitriangular structures) obtained by Belavin and Drinfeld, and reduces to it for $\l=0$ (see \cite{BD}).\\
2. When $\l=\h$ one recovers the classification result Theorem 4.2: the only $\h$-graded generalized
Belavin-Drinfeld triples are of the form $\Gamma_1=\Gamma_2, \tau=Id$, and in this case the r-matrix
(\ref{E:app}) corresponds to $r_X^1(\lambda)$, with $X=\Gamma_1$.\\ 3. Theorem~\ref{T:app} is proved in
\cite{S} under the additional assumption that $\l$ contains a regular semisimple element.
However, the proof easily extends to the present situation.

\section{Appendix B: Proof of the ABRR relation for $U_q(\g)$} 
\paragraph{}We keep the notations of Section 8. Recall the Drinfeld construction of 
the quantum Casimir element of $U_q(\g)$. 
Let ${\mathcal R}$ be the universal R-matrix for $U_q(\g)$.
Let us write ${\mathcal R}=\sum_i a_i \otimes b_i$ and set $u=\sum S(b_i)a_i$. Then $u=q^{2\rho}z$ where $z$ is a central element
in a completion of $U_q(\g)$, which is called the quantum Casimir
element. Moreover, for any $\mu \in \h^*$ we have 
$$ux_\mu=q^{-\sum_i x_i^2}x_\mu=q^{-(\mu,\mu)}x_\mu$$
hence $zx_\mu=q^{-(\mu,\mu+2\rho)}x_\mu$ and $u_{|M_\mu}=q^{-(\mu,\mu+2\rho)}q^{2\rho}$.\\
\hbox to1em{\hfill}Now let $V$ and $W$ be two finite-dimensional $U_q(\g)$-modules, $v \in V$ and $w \in W$ homogeneous vectors of weight $\mu_v$ and $\mu_w$ respectively, and consider the expectation value
$$X_{vw}(\mu)=\langle x^*_{\mu-\mu_v-\mu_w},\Phi^w_{\mu-\mu_v}u_{|M_{\mu-\mu_v}}\Phi^v_\mu x_\mu \rangle.$$
We will compute $X_{vw}(\mu)$ in two different ways. On one hand, 
\begin{equation*}
\begin{split}
X_{vw}(\mu)&=\langle x^*_{\mu-\mu_v-\mu_w},\Phi^w_{\mu-\mu_v}(q^{2\rho}z)_{|M_{\mu-\mu_v}}\Phi^v_\mu x_\mu \rangle\\
&=\langle  x^*_{\mu-\mu_v-\mu_w},\Phi^w_{\mu-\mu_v} z_{|M_{\mu-\mu_v}}q^{-2\rho}_{|V}\Phi^v_\mu q^{2\rho}x_\mu \rangle\\
&=q^{-(\mu-\mu_v,\mu-\mu_v+2\rho)}q_{|V}^{2\rho}q^{(2\rho,\mu)}J(\mu)(w \otimes v)\\
&=q^{-(\mu-\mu_v,\mu-\mu_v)+2(\rho,\mu_v)}q_{|V}^{-2\rho}J(\mu)(w \otimes v)
\end{split}
\end{equation*}
On the other hand, we have $(1 \otimes \Delta^{op}) {\mathcal R}={\mathcal R}^{12}{\mathcal R}^{13}$ hence, by the intertwining property,
$$(1 \otimes \Delta) (1 \otimes S) {\mathcal R}=\sum_{ij} a_ia_j \otimes S(b_i)\otimes S(b_j).$$
Thus
\begin{equation}\label{E:app1}
X_{vw}(\mu)=\sum_{i,j} S(b_j)_{|W} \langle x^*_{\mu-\mu_v-\mu_w}, S(b_i)_{|M_{\mu-\mu_v-\mu_w}} \Phi^w_{\mu-\mu_v} a_ia_{j|M_{\mu-\mu_v}} \Phi^v_\mu x_\mu \rangle.
\end{equation}
Now, since $x^*_{\mu-\mu_v-\mu_w}$ is a lowest weight vector and since $${\mathcal R}\in 
(1+U'_q(\b_-) \otimes U'_q(\b_+)) q^{\sum_i x_i \otimes x_i}$$
equation (\ref{E:app1}) reduces to
\begin{equation}\label{E:app2}
\begin{split}
X_{vw}(\mu)&=\sum_{j} S(b_j)_{|W} \langle x^*_{\mu-\mu_v-\mu_w}, q^{-\sum x_i^2}_{|M_{\mu-\mu_v-\mu_w}}q^{-\sum x_i \otimes x_i}_{|M_{\mu-\mu_v-\mu_w}\otimes W} \Phi^w_{\mu-\mu_v} a_{j|M_{\mu-\mu_v}} \Phi^v_\mu x_\mu \rangle\\
&=q^{-(\mu-\mu_v-\mu_w)^2}\sum_j \big(S(b_j)q^{-\mu+\mu_v+\mu_w}\big)_{|W}\langle x^*_{\mu-\mu_v-\mu_w} \Phi^w_{\mu-\mu_v} a_{j|M_{\mu-\mu_v}} \Phi^v_\mu x_\mu \rangle
\end{split}
\end{equation}
Using the relation
$$(\Delta \otimes S){\mathcal R}=\sum a_l \otimes a_k \otimes S(b_k)S(b_l)$$
we have
$$\sum_j S(b_j)_{|W} \otimes \Phi^v_\mu
a_{j|M_\mu}=\bigg(\sum_{k,l}
\big(S(b_k)S(b_l)\big)_{|W} \otimes a_{l|M_{\mu-\mu_v}} \otimes
a_{k|V} \bigg)\Phi^v_\mu,$$ i.e
$$\big(\sum S(b_l)_{|W} \otimes a_{l|M_{\mu-\mu_v}}\big)\Phi^v_\mu=\big(\sum S(b_k)_{|W} \otimes a_{k|V}\big)^{-1} \sum S(b_j)_{|W} \otimes \Phi^v_\mu a_{j|M_\mu}.$$
Substitution of this in (\ref{E:app2}) gives
\begin{equation}\label{E:app3}
X_{vw}(\mu)=q^{-(\mu-\mu_v-\mu_w)^2}\big(\sum S(b_k)_{|W} \otimes a_{k|V}\big)^{-1}q_{|W}^{-2\mu+\mu_v+\mu_w} J(\mu)(w \otimes v)
\end{equation}
\paragraph{Claim:} we have $\big(\sum S(b_k) \otimes a_{k}\big)^{-1}=(1 \otimes q^{2\rho}) {\mathcal R} (1 \otimes q^{-2\rho})$.\\
\textit{Proof.} We have $S^2(a)=q^{2\rho}a q^{-2\rho}$ for any $a
\in U_q(\g)$. Hence $\sum a_k \otimes S(b_k)=(1 \otimes
q^{2\rho})(1 \otimes S^{-1}) (1 \otimes q^{-2\rho})$. The claim now
follows from the relation\\ $(1 \otimes S^{-1}){\mathcal
R}={\mathcal R}^{-1}$.
\paragraph{}Thus, combining (\ref{E:app1}) and (\ref{E:app3}), we get
\begin{equation*}
\begin{split}
q^{-(\mu-\mu_v,\mu-\mu_v+2\rho)}&q^{-2\rho}_{|V}
q^{2(\mu,\rho)}J(\mu)(w \otimes
v)\\
=&q^{-(\mu-\mu_v-\mu_w)^2}q^{2\rho}_{|W}{\mathcal
R}^{21}q_{|W}^{\mu_v+\mu_w-2(\mu+\rho)}J(\mu)(w\otimes v).
\end{split}
\end{equation*}
This implies that
\begin{equation}\label{E:app35}
\begin{split}
q^{-(\mu_v+\mu_w)^2}&{\mathcal
R}^{21}q^{\mu_v+\mu_w-2(\mu+\rho)}_{|W}J(\mu)(w \otimes v)\\
=&q^{\mu^2-2(\mu,\mu_v+\mu_w)}q^{-2\rho}_{|V}q^{-2\rho}_{|W} q^{2(\mu,\rho)}q^{-(\mu-\mu_v,\mu-\mu_v+2\rho)}J(\mu)(w \otimes v)
\end{split}
\end{equation}
Using the weight zero property of ${\mathcal R}$, we can rewrite the l.h.s of this last equation as
\begin{equation}\label{E:app4}
{\mathcal R}^{21}q^{-\sum x_i \otimes x_i}_{|W \otimes V} q^{-\sum x_i^2}_{|V} q^{-2(\mu+\rho)}_{|W}J(\mu)(w \otimes v)
\end{equation}
Similarly, using the weight zero property of $J(\mu)$, it is easy to see that the right hand side is equal to
\begin{equation}\label{E:app5}
q^{-2(\mu+\rho,\mu_w)}q^{-\mu^2}J(\mu)=J(\mu)q^{-2(\mu+\rho)}_{|W}q^{-\sum x_i^2}_{|V}(w \otimes v)
\end{equation}
The ABRR equation now follows from (\ref{E:app35}), (\ref{E:app4}), (\ref{E:app5}) and the weight zero property of $J(\mu)$.

\section{Appendix C: Extension of dynamical 2-cocycles}
\paragraph{}Let $H$ be a quantized enveloping algebra, i.e $H$ is a $\C[[\hbar]]$-Hopf algebra which is topologically free as a $\C[[\hbar]]$-module and such
that $H/\hbar H \simeq U\g$ for some Lie algebra $\g$. Hence $H \simeq U\g[[\hbar]]$ with deformed algebra and coalgebra structures. Let $\h \subset \g$ and suppose that the embedding $i:U(\h)[[\hbar]] \to H$ is a Hopf algebra map. In
particular, $\h \subset H$ consists of primitive elements.
\paragraph{}Let $J(\lambda):\;\h^* \to (H {\otimes} H)^\h$ be a meromorphic
solution of the dynamical 2-cocycle equation :
\begin{equation}\label{E:C1}
J^{12,3}(\lambda)J^{12}(\lambda-\hbar h^{(3)})=J^{1,23}(\lambda)J^{23}(\lambda).
\end{equation}
Here ${\otimes}$ denotes the tensor product in the category 
topologically free $\C[[\hbar]]$-modules, and 
$J^{12}(\lambda-\hbar h^{(3)})$ is given by the Taylor expansion: 
$$ J^{12}(\lambda-\hbar h^{(3)})\stackrel{def}{=} J(\lambda) - \hbar
\sum_i \frac{\partial J}{\partial y_i}\big(\lambda\big) x_i +
\ldots,$$ where
$(y_i)$ is a basis of $\h^*$ and $(x_i)$ the dual basis of $\h$.

\paragraph{}
It is easy to see that if 
f $J$ satisfies the dyanmical 2-cocycle equation, 
then 
$J(\lambda-\mu)$ can be explicitly expressed in
terms of $J(\lambda)$ for any $\mu$.
Namely, let us extend any linear function $\mu$ on $\h$ to a ring
homomorphism $U(\h) \to \C[[\hbar]]$. 
Then $J(\lambda-\mu)$ is determined through the formula
$$J^{12}(\lambda-\mu)=\big(\frac{\mu}{\hbar}\big)_3\big(J^{12,3}(\lambda)^{-1}J^{1,23}(\lambda)J^{23}(\lambda)\big)$$
This means that $J(\lambda)$
is completely determined by its value at any point. 

The main result of this subsection is the following Proposition. 

Let $U(\hbar \h)$ denote the enveloping algebra of $\hbar\h$,
i,e the $\C[[\hbar]]$-subalgebra of $H$ spanned by elements $\hbar h$,
$h \in \h$. 
\paragraph{}

\begin{prop}Let $J \in (H {\otimes} H)^\h$ satisfy the relation
\begin{equation}\label{E:C2}
(J^{12,3})^{-1}J^{1,23}J^{23} \in H {\otimes} H {\otimes} U(\hbar\h).
\end{equation}
Set 
$$J(-\mu)=\frac{\mu_3}{\hbar}\big((J^{12,3})^{-1}J^{1,23}J^{23}\big)\;\qquad\;for\;all
\;\mu\in \h^*.$$
Then $J(\lambda)$ satisfies the dynamical 2-cocycle condition (\ref{E:C1}).
This establishes a one-to-one correspondence between solutions
of (\ref{E:C1}) in $(H \otimes H)^\h$ regular at $0$ and solutions of (\ref{E:C2}) in $(H \otimes H)^\h$.
\end{prop}
\noindent

Note that the right hand side 
of the definition of $J(-\mu)$ makes sense since
the action of
$\frac{\mu}{\hbar}$ on $U(\hbar \h)$ is well-defined.

\textit{Proof.} The proof is a direct computation. By definition we have
\begin{align*}
J^{23}(\lambda)&=\big(-\frac{\lambda}{\hbar}\big)_4\big((J^{23,4})^{-1}J^{2,34}J^{34}\big)\\
J^{1,23}(\lambda)&=\big(-\frac{\lambda}{\hbar}\big)_4\big((J^{123,4})^{-1}J^{1,234}J^{23,4}\big)\\
(J^{12,3}(\lambda))^{-1}&=\big(-\frac{\lambda}{\hbar}\big)_4\big((J^{34})^{-1}(J^{12,34})^{-1}
J^{123,4}\big).
\end{align*}
Thus
$$\big(\frac{\mu}{\hbar}\big)_3\big((J^{12,3})^{-1}(\lambda)J^{1,23}(\lambda)J^{23}(\lambda)\big)
=\big(\frac{\mu}{\hbar}\big)_3\big(-\frac{\lambda}{\hbar}\big)_4\big((J^{34})^{-1}(J^{12,34})^{-1}J^{1,234}J^{2,34}J^{34}
\big).$$
But
$$(J^{12,34})^{-1}J^{1,234}J^{2,34}=\Delta_3 \big((J^{12,3})^{-1}J^{1,23}
J^{23}\big) \in H {\otimes} H {\otimes} \Delta(U(\hbar\h))$$
and $J \in (H {\otimes} H)^{\h}$. Hence
$$(J^{34})^{-1}(J^{12,34})^{-1}J^{1,234}J^{2,34}J^{34}=(J^{12,34})^{-1}
J^{1,234}J^{2,34}=\Delta_3 \big((J^{12,3})^{-1}J^{1,23}J^{23}\big).$$
Moreover, for any $h \in \h$, $\big(\frac{\mu}{\hbar}\big)_3\big(-\frac{\lambda}{\hbar}\big)_4\Delta(h)=\big(\frac{\mu}{\hbar}\big)_3\big(-\frac{\lambda}{\hbar}\big)_4
(h \otimes 1 + 1 \otimes h)=\big(\frac{\mu-\lambda}{\hbar}\big)h$. Therefore
$$\big(\frac{\mu}{\hbar}\big)_3\big((J^{12,3}(\lambda))^{-1}J^{1,23}(\lambda)J^{23}(\lambda)\big)
=\big(\frac{\lambda-\mu}{\hbar}\big)_3\big((J^{12,3})^{-1}J^{1,23}J^{23}\big)=J^{12}(\lambda-\mu),$$
which proves the Proposition. \qed

\section{Appendix D. Fusion matrices and Shapovalov forms}
\paragraph{}Let $\g$ be a semisimple complex Lie algebra. Let $\h \subset \g$
be a Cartan subalgebra. Consider the Chevalley involution $\omega$ defined by
$$\omega(e_\alpha)=-e_{-\alpha},\qquad \omega(e_{-\alpha})=-e_{\alpha},\qquad
\omega(h_\alpha)=-h_{\alpha}$$
for every simple root $\alpha$. Let $U$ stand for
either of the algebras $U(\g)$ or $U_q(\g)$, and let $U^+=U(\n_+)$ if $U=U(\g)$ and $U^+=U_q(\n_+)$ if $U=U_q(\g)$. It is clear that $\omega$ extends
to an involutive automorphism of $U$.
\paragraph{}Let $\lambda \in \h^*$ and let $M_\lambda^+$ be the associated
highest weight Verma module. It is well-known that there exists a unique
symmetric bilinear form $\langle\,,\,\rangle$ on $M^+_\lambda$ such that
$$\langle x^+_\lambda,x^+_\lambda \rangle=1,\qquad \langle e_\alpha u,v\rangle=
\langle u,e_{-\alpha}v\rangle,\qquad \forall\;u,v\in M^+_\lambda.$$
This form is called the \textit{Shapovalov form}. Furthermore, $\langle\,,\,
\rangle$ is nondegenerate if and only if $M^+_\lambda$ is irreducible. In
particular, it is nondegenerate for generic values of $\lambda$.\\
\hbox to1em{\hfill}From now on, we will assume that $\lambda$ is generic. Then
we have a $U$-module isomorphism $(M_\lambda^+)^*\stackrel{\sim}{\to}
M_{-\lambda}^-$ where $M^-_{-\lambda}$ is the lowest weight Verma module
of lowest weight $-\lambda$ and $(M_\lambda^+)^*$ is the restricted dual of
$M_\lambda^+$. Furthermore, $(M^-_{-\lambda})^\omega=M^+_\lambda$
where $\omega$ denotes the operation of twisting by $\omega$, and the natural
pairing $M^+_{\lambda} \otimes M^-_{-\lambda} \to \C$ is identified with the
Shapovalov form. Let $X_\lambda \in M^+_{\lambda} \hat{\otimes} M^-_{-\lambda}$
be the inverse element to this form. 
\paragraph{}The following result relates $X_\lambda$ with the fusion matrix.
\begin{prop} Let $J(\lambda)$ be the fusion matrix. Then
$X_\lambda=J(0)(x^+_\lambda \otimes x^-_{-\lambda})$.
\end{prop}
\noindent
\textit{Proof.} Consider the composition of intertwiners
$$M_0^+ \stackrel{\Phi^{x_\lambda^+}_0}{\longrightarrow} M^+_\lambda \hat{\otimes}
M^-_{-\lambda}\stackrel{\Phi^{x_{-\lambda}^-}_\lambda}{\longrightarrow}
\big(M_0^+ {\otimes} M_\lambda^+\big)\hat{\otimes}
M^-_{-\lambda}.$$
Note that the first intertwiner is well defined because for 
generic $\lambda$ the module $M_\lambda^+$ is irreducible, 
and the second one obviously exists for any $\lambda$. 
The expectation value $\langle x_0^{+*},\Phi^{x_\lambda^+}_0
\Phi^{x_{-\lambda}^-}_\lambda x^+_0\rangle$ is by definition
$J(0)(x_\lambda^+ \otimes x^-_{-\lambda})$. On the other hand, note
that $\Phi^{x_\lambda^+}_0(x_0^+)=X_\lambda$. Indeed,
$\Phi^{x_\lambda^+}_0(x_0^+)$ is the unique $U^+$-invariant element of
$ M^+_\lambda \hat{\otimes}M^-_{-\lambda}$ of the
form $x_\lambda^+ \otimes x_{-\lambda}^- +\;l.o.t$, and $X_\lambda$ satisfies
both conditions. Thus  $\langle x_0^{+*},\Phi^{x_\lambda^+}_0
\Phi^{x_{-\lambda}^-}_\lambda x^+_0\rangle=X_\lambda$. \qed

\paragraph{Remark.} The above proposition admits a direct
generalization to the Kac-Moody setting.
\paragraph{Remark.} Consider the subalgebra $\widetilde{U_q({\n}_-)}$ generated by the left components
of $\mathcal{R}_0=\mathcal{R}q^{-\sum_i x_i \otimes x_i}$ (see
Section 8.1).
It follows from the ABRR equation that 
$J(\lambda) \in \widetilde{U_q({\n}_-)} \hat{\otimes}
U_q(\b_+)$. This implies, in light of Appendix C, that 
Proposition 13.1 gives a complete characterization of the fusion
matrix. Thus, the problems of computing the Shapovalov form and fusion matrix
are equivalent.

\section{Review of literature}

In conclusion, we would like to give a brief review 
of the existing literature on the dynamical Yang-Baxter equations. 
We would like to make it clear that this list is by no means complete, 
and contains only some of the basic references which are relevant to 
this paper.  
 
The physical paper in which the dynamical Yang-Baxter equation 
was first considered is \cite{GN}; dynamical R-matrices are also discussed in 
\cite{Fad1},\cite{AF}. 

The classical dynamical Yang-Baxter equation and examples of its
solutions were introduced in \cite{F}. Its geometric interpretation in
terms  of Poisson groupoids of Weinstein \cite{W} was introduced in
\cite{EV1}.  Solutions of this equation were studied and classified in
\cite{EV1},
\cite{S}. The relationship of solutions of this equation to Poisson groupoids 
and Lie bialgebroids was further explored in \cite{LX} and \cite{BKS}.
The relationship of solutions of the classical dynamical Yang-Baxter equation 
(defined on noncommutative Lie algebras) to equivariant cohomology 
is discussed in \cite{AM}. The relationship to integrable systems is 
discussed in \cite{ABB}.

Quantum groups associated to a dynamical R-matrix
were first introduced in \cite{F}. 
In the case when the R-matrix is elliptic, 
they are called elliptic quantum groups.
These quantum groups and their representation theory 
(for the Lie algebras $sl_n$), as well as their relationship 
with integrable systems, were 
systematically studied in 
\cite{FV1,FV2,FV3}. The papers \cite{EV2,EV3} study 
the trigonometric versions of these quantum groups (for any 
simple Lie algebra). 

Quantum groupoids were introduced by Maltsiniotis (in the case
when the base 
is classical), and by Lu in \cite{Lu}  in the general case. 
The interpretation of dynamical quantum groups as quantum groupoids 
was first discussed in \cite{EV2,EV3}, and further studied in 
\cite{Xu1}, \cite{Xu2}. 
The interpretation of dynamical quantum groups as 
quasi-Hopf algebras is discussed in \cite{BBB}, and was further developed in 
\cite{Fr,JKOS}, \cite{ABRR},\cite{EF}. 
The connection between these 
two interpretation is discussed in \cite{Xu1}.
The theory of dynamical quantum groups at roots of 1 and their
interpretation via 
weak Hopf algebras was discussed in 
\cite{EN}. 

Quantum KZB equations (which are not discussed in these notes) 
were introduced in \cite{F}, and studied in \cite{FTV1,FTV2,MV,FV3,FV4,FV5,
FV6,FV7}. 
Monodromy of quantum KZ equations \cite{FR}, which 
yields dynamical R-matrices of the elliptic quantum groups, 
is computed in \cite{TV1,TV2}. 

The theory of traces of intertwining operators 
for Lie algebras and quantum groups and its applications 
to the theory of special functions (in particular, Macdonald theory) 
is developed in 
\cite{B},\cite{F},\cite{E},\cite{EK1,EK2,EK3,EK4},\cite{K1,K2,K3}, 
\cite{EFK}, \cite{ES1,ES2}.
The relationship of this theory with dynamical 
classical r-matrices and quantum R-matrices 
is studied in \cite{EV4}, \cite{ESch1},\cite{ESch2}. 

\small{
\end{document}